\documentclass[10pt]{elsarticle}

\usepackage{amsmath,amssymb,amsfonts,color,amsthm}
\usepackage{mathtools}
\usepackage{algorithmic} 
\usepackage{algorithm,booktabs}
\usepackage{tikz}
\usepackage[mathscr]{eucal}

\newtheorem{theorem}{Theorem}[section]
\newtheorem{lemma}[theorem]{Lemma}

\newtheorem{proposition}[theorem]{Proposition}

\newtheorem{definition}[theorem]{Definition}

\newtheorem{rmk}{Remark}[section]

\newcommand{\R}{\mathbb{R}}
\newcommand{\N}{\mathbb{N}}

\newcommand{\red}[1]{\textcolor{black}{#1}}
\newcommand{\blue}[1]{\textcolor{black}{#1}}
\newcommand{\new}[1]{\textcolor{black}{#1}}

\newcommand*{\Cdot}[1][1.25]{%
  \mathpalette{\CdotAux{#1}}\cdot%
}
\newdimen\CdotAxis
\newcommand*{\CdotAux}[3]{%
  {%
    \settoheight\CdotAxis{$#2\vcenter{}$}%
    \sbox0{%
      \raisebox\CdotAxis{%
        \scalebox{#1}{%
          \raisebox{-\CdotAxis}{%
            $\mathsurround=0pt #2#3$%
          }%
        }%
      }%
    }%
    \dp0=0pt %
    \sbox2{$#2\bullet$}%
    \ifdim\ht2<\ht0 %
      \ht0=\ht2 %
    \fi
    \sbox2{$\mathsurround=0pt #2#3$}%
    \hbox to \wd2{\hss\usebox{0}\hss}%
  }%
}

\newcommand*{\dummy}{\,\Cdot[1.5]\,}

\journal{Advances in Computational Mathematics}

\begin{document}

\begin{frontmatter}



%via\\ Proper Orthogonal Decomposition and \\Dynamic Mode Decomposition}


\title{A Certified Model Reduction Approach for Robust Parameter Optimization with PDE Constraints }

\author{Alessandro Alla}
\ead{alla@mat.puc-rio.br}
\address{Department of Mathematics, PUC-Rio, R. Mq. S. Vicente 225, Rio de Janeiro 22453-900, Brazil.}

\author{Michael Hinze}
\ead{michael.hinze@uni-hamburg.de}
\address{Department of Mathematics,
Universit\"at Hamburg, Bundesstr. 55, 20146 Hamburg, Germany.}

\author{\new{Philip Kolvenbach}}
\ead{kolvenbach@mathematik.tu-darmstadt.de}
\address{\new{Department of Mathematics, Technische Universit\"at Darmstadt,
 Dolivostr. 15, 64293 Darmstadt, Germany.}}

\author{Oliver Lass}
\ead{lass@mathematik.tu-darmstadt.de}
\address{Department of Mathematics, Technische Universit\"at Darmstadt,
 Dolivostr. 15, 64293 Darmstadt, Germany.}

\author{Stefan Ulbrich}
\ead{ulbrich@mathematik.tu-darmstadt.de}
\address{Department of Mathematics, Technische Universit\"at Darmstadt,
 Dolivostr. 15, 64293 Darmstadt, Germany.}

\begin{abstract}
We investigate an optimization problem governed by an elliptic partial differential equation with uncertain parameters. We introduce a robust optimization framework that accounts for uncertain model parameters. The resulting non-linear optimization problem has a bi-level structure due to the min-max formulation. To approximate the worst-case in the optimization problem we propose linear and quadratic approximations. However, this approach still turns out to be very expensive, therefore we propose an adaptive model order reduction technique which avoids long offline stages and provides a certified reduced order surrogate model for the parametrized PDE which is then utilized in the numerical optimization. Numerical results are presented to validate the presented approach.

 \end{abstract}


\end{frontmatter}

\section{Introduction}
Parameter optimization governed by partial differential equations (PDEs) is a well-studied topic due to its relevance in industrial applications. If the model is considered to be {\em perfect} one can perform standard techniques to solve these type of problems (see e.g. \cite{HPUU09,Tro10} and the references therein). However, its numerical approximation can be very expensive due to the dimension of the discretized PDE. For this reason, in the last decade, model order reduction techniques were introduced and successfully applied in the context of PDE constrained optimization. Model order reduction works in a Galerkin projection framework, where the basis functions are non-local and built upon information of the underlying system. Although a detailed description of these contributions goes beyond the scope of this paper, we want to mention  Proper Orthogonal Decomposition (POD) and the reduced basis (RB) method.
The former works in general situations such as time-dependent problems, and parametric steady and unsteady equations (see e.g. \cite{V13} for a presentation of the method within different applications). The latter is mainly used in the context of e.g. many-query scenarios (see e.g. \cite{RHP08}) for parametric steady problems where the basis functions are selected by means of a greedy algorithm. It is also possible to combine RB and POD in the so-called POD-Greedy algorithm \cite{HO08} for parametric unsteady models. The strength of these methods is the presence of an a-posteriori error estimator which certifies the quality of the surrogate models. Model order reduction has been applied successfully to parameter optimization with PDE constrains in e.g. \cite{DH15,OP07}. However, we note that all these works strongly rely on the offline-online decomposition. In the current work, we propose to build the surrogate model towards the optimization. Similar approaches were studied in e.g. \cite{QGVW16,ZF14}, \blue{but our approach does not require offline-online decomposition} and the reduced order model gets updated to desired accuracy during the optimization process.

%{Model order reduction has been applied successfully to optimal control problems in both open-loop (e.g. \cite{HV05,KG14,SV10}) and closed-loop (e.g. \cite{AV15,AK01,KVX04}) frameworks.}

For the purpose of this work we assume that the PDE which governs our system is given but material imperfections are present. These are due to e.g. manufacturing and introduce uncertainty to the model. This problem arises in many real life applications. One way to include this uncertainty into the optimization process is through robust optimization. In this case no probabilistic model of the uncertainty is required. Instead, a deterministic approach is applied by assuming that the uncertainty is restricted to a bounded uncertainty set. Using the notion of a robust counterpart the associated original uncertain optimization problem is reformulated. The solution obtained in this way stays feasible for all realizations from the uncertainty set and, at the same time, yields the best worst-case value of the objective function. Techniques which resolve the uncertainty by means of stochastic optimization e.g. \cite{BL97,KS16,Oks05,TKXP12} are an alternative to the presented approach. These methods depend on sampling the uncertain parameters and hence can become prohibitive expensive in the context of PDE constrained optimization. For a general discussion of robust optimization we refer to e.g. \cite{BGN09,BN02,BBC11}.

The idea is to utilize a suitable approximation of the robust counterpart, e.g. \cite{DBK06,HD13,LU17,Sic13,Zha07} or exploit specific properties e.g. \cite{HD13}. We investigate the approximation of the robust counterpart using a quadratic model. This allows us to reformulate the robust optimization problem as a {\it mathematical program with equilibrium constraints} (MPEC). This approach has been investigated in \cite{LU17,Sic13, AU16} in the context of PDE constrained optimization problems, where \cite{LU17} forms the starting point for the present investigations.
Our model will be a linear elliptic parametric equation where an affine decomposition will be applied to work with a reference parameter. The worst-case problem leads a non-linear optimization problem with a $\min-\max$ formulation. After approximating the inner maximization derivative based optimization is applied. In the setting of our work we will utilize a sensitivity based approach since we assume to deal with only a
few parameters. This approach is computationally expensive because it requires the solution of multiple PDEs. 

In this work the focus is on solving the robust counterpart in an efficient way utilizing a POD-based reduced order method. More specifically, we propose a certified robust optimization procedure which, following an idea proposed in \cite{AH01}, successively adapts the POD model during the iterative adjustment of the optimal robust design through appropriate enrichment of the snapshot set and avoids expensive offline stages.  The method is certified by an a-posteriori error estimator for the state variable and the sensitivities. Therefore, error estimators for the state and the sensitivities are required. A generalized error estimator is derived to cover these needs.
Here we build upon the robust optimization framework developed in \cite{LU17}, which for the convenience of the reader is summarized in Section~\ref{sec:opt}\new{, and which we have now complemented with a moving expansion point approach to further improve the accuracy of the utilized approximations, see Section 3.2.3 for details.}
%\new{We extend this framework with a new method to increase the accuracy of the utilized approximations.}
Our approach \new{to model reduction} complements the method proposed in \cite{LU17}. In fact, in
\cite{LU17} the model reduction in the application is performed by a POD-greedy
procedure to reduce the degrees of freedom related to the rotation angle of
a \new{quasistatic} model of an electrical machine. This 
POD surrogate model is then used to speed up the robust optimization.
In this paper, we generate the reduced order model by applying POD to a
snapshot set resulting from state and sensitivities corresponding to previous
optimization iterates. The quality of the obtained optimal solution is
evaluated by a posteriori error estimators. As long as the error is above a
specified level, the snapshot set is updated by the state and sensitivities corresponding
to the current iterate of the optimization variables.
We demonstrate the performance and efficiency of our
approach for robust design of electrical machines, where we use the setting
of \cite{LU17} in the static case, i.e.,  we do not consider rotation.

The paper is organized as follows: in Section \ref{sec2} we present the mathematical model and in Section \ref{sec:opt} we present the nominal and the robust optimization problem. Section \ref{sec:pmor} focuses on model order reduction for the optimization problem and in Section \ref{sec:test} we illustrate the effectiveness of the discussed methods by numerical examples.

% -------------------------------------------------
% -------------------------------------------------
% --------------- SECTION 2 -----------------
% -------------------------------------------------
% ------------------------------------------------- 
\section{Elliptic Parametric PDEs}\label{sec2}
We deal with an abstract linear stationary
equation governed by a uniformly continuous and uniformly coercive bilinear
form $\tilde a$ stemming e.g.  from a parametric elliptic PDE.  The
parameter sets are denoted by $\mathcal D \subset \R^{N_p}$ and $\mathcal
U_k \subset \R^{N_\Phi}$, where parameters $p\in \mathcal D$ enter the
problem via a $p$-dependent regular bounded domain $\Omega(p) \subset \R^2$
which give rise to a $p$-dependent real and separable Hilbert space
$X(\Omega(p))$ and dual space $X'(\Omega(p))$.  The parameters $\Phi \in
\mathcal U_k$ will model the uncertainty in the problem.  We note that for
our abstract setting we impose the by now standard assumptions on parameter
separability which are formulated for example in \cite{H17}. Those assumptions are in
particular met for the application considered in Section \ref{sec:test}. 
Our abstract problem then reads

%We note that we follow the assumptions provided in \cite{H17} which later allow the parameter separability and the use of the $\min-\hat\Phi$ theorem. Our abstract problem then reads}
%
%\myc{
%ASSUMPTIONS FROM \cite{H17}
%\begin{enumerate}
%\item $X$ a real, seprable Hilbert space with inner product $\langle\cdot,\cdot\rangle$, norm $\|\cdot\|$ and dual space $X'$ with norm $\|\cdot\|_{X'}$
%\item parameter domain $\mathcal{P}\subset\R^p$
%\item a parameter dependent bilinear form $a(\cdot,\cdot;\mu)$ and linear form $f(\cdot,\mu)$ for all $\mu\in\mathcal{P}.$
%\item We assume the bilinear form and the linear forms to be uniformly continuous and the bilinear form to be uniformly coercive 
%\end{enumerate}
%}
%
%
%\begin{figure}[htbp]
%\includegraphics[scale=0.5]{parameter_separability}
%\end{figure}
%In this section we introduce our abstract model problem: we deal with a linear elliptic equation with arbitrary \blue{well-defined} boundary conditions, considering $\mathcal{D}\subset \R^{\red{N_p}},\,\red{N_p}\ge1$ as the parameter space and $\Omega(p)\blue{\subset\R^2}$ as a domain:
%
\begin{equation}\label{pde_par}
\left\{\begin{array}{l}
\hbox{For }p\in\mathcal{D} \hbox{ find }\tilde{u}\in X(\Omega(p))\hbox{ s.t}\\
\tilde{a}(\tilde{u},v;p)=\tilde{f}(v;p,\phi),\quad \forall v\in X(\Omega(p)),
\end{array}\right.
\end{equation}
where $\tilde{u}$ is the unknown variable.
% $X$ is a suitable Hilbert Space, $X'$ its dual space, $p\in \mathbb R^{N_p}$ a $N_p$ dimensional parameter, \blue{$\Omega(p)$} a regular bounded parameter dependent domain, , \blue{$\phi\in\mathcal{U}_k$ are parameters, with $\mathcal{U}_k\subset \mathbf R^{N_\phi}$ which later will model the uncertainty in our mathematical model.}. 
In our application uncertainty only enters through the right hand side,
where $\tilde a$ and $\tilde{f}$ are given by %the bilinear bounded and coercive form $\tilde{a}$ and right hand side $\tilde{f}$ are given by
$$
\tilde{a}(w,v;p) = \int_{\Omega(p)} \nu \nabla w \cdot \nabla v \,\mathrm dx
\quad\mbox{and}\quad
\tilde{f}(v;p,\phi) = \int_{\Omega(p)} \blue{\tilde{f}(\phi)} v \,\mathrm dx,
$$
where $\nu\ge \bar{\nu}>0$ is an isotropic material coefficient (see Section \ref{sec:opt} for more details.). %, i.e., $\nu(\phi) = \mathrm{diag}(\nu^1(\phi),\nu^1(\phi))$ and $\nu^1(\phi)>0$.}
 Note that under these assumptions problem \eqref{pde_par} admits a unique solution.
%\myc{MH2: Furthermore $\nu(\phi)$ this indicates that only this constant later will be uncertain. Is it this what we mean? }

The function space $X(\Omega(p))$  is such that $H^1_0(\Omega(p))\subset \blue{X(\Omega(p))} \subset H^1(\Omega(p)),$ with %$H^1(\Omega(p))$ the Sobolev Space (please refer to \cite{Eva08} for more details on partial differential equations) defined as:
$$
H^1(\Omega(p)):=\left\{f\in L^2(\Omega(p)): D^\alpha f\in (L^2(\Omega(p)))^2, |\alpha| \leq1\right\},
$$
where $f$ is a measurable function, \blue{$\alpha\in\mathbb{N}_0^2$}, $D^\alpha f$ denotes the weak $\alpha^{th}-$ partial derivative of $f$ and the functional spaces $L^2(\Omega(p))$, $H^1_0(\Omega(p))$ are defined as follows:
$$
L^2(\Omega(p)):=\left\{f:\Omega(p)\rightarrow\R, \int_{\Omega(p)} f(x)^2\, \mathrm dx<\infty\right\}, 
$$
$$
H^1_0(\Omega(p)):=\left\{f\in H^1(\Omega(p)): f\equiv0 \hbox{ a.e.  on }\partial\Omega(p) \hbox{ in the sense of traces} \right\}.
$$
In order to obtain a computationally fast model and to avoid re-meshing when the parameter changes we require parameter separability,
see, e.g., \cite{H17}, and assume that the domains $\Omega(p)$ in \eqref{pde_par} can be transformed to a fixed domain $\Omega(\bar p):= \cup_{q=1}^ Q \Omega_q(\bar p)$, \blue{where $\bar p$ denotes the reference parameter.}
Without loss of generality, we assume that the domain of interest $\Omega(p)$ can be decomposed in $Q$ non-overlapping triangles and the transformation on each triangle \red{is linear, whereas piecewise-linear} and continuous over the whole domain  according to:
\begin{align}\label{exp_dec}
\mathcal{T}_q( x; p): \Omega_q(\bar p) &\rightarrow \Omega(p)\nonumber\\
  x &\blue{{}\mapsto{}} C_q(p) x + d_q(p), % ATTENTION: delete the two pairs of braces {} when deleting \blue
\end{align}
for $q=1,\ldots,Q,$ where $C_q(p)\in\R^{2\times 2}$ and $d_q(p)\in\R^2$. 
\blue{From here onwards we write $\bar \Omega$ instead $\Omega(\bar p)$ whenever this is appropriate.}%\blue{We also note that whenever we want to ease the notation we will use $\bar\Omega$ for $\Omega(\bar p)$.}
% $\mathcal T(p): \bar \Omega \rightarrow \Omega(p)$ such that $\Omega(p) = \mathcal T(p)\bar\Omega.$}
%which transforms problem \eqref{pde_par} to a fixed reference domain $\bar \Omega$ by means of a \myc{MH:piecewise affine and continuous?} linear transformation. For this purpose, let us assume that there is a piecewise-affine transformation
%
%$$
%\mathcal{T}:\bar{\Omega}\rightarrow\Omega(p),
%$$
%
%such that $\Omega(p)=\mathcal{T}(p)\bar{\Omega}$. 

As it is
shown, e.g., in \cite{RHP08}, the linear parametric map $\mathcal{T}_q(p)$ and its Jacobian $J_{\mathcal{T}_q}$ allow the definition of the bilinear and linear forms on the reference domain $\bar{\Omega}$. Then, problem \eqref{pde_par} can be rewritten as
\begin{equation}\label{pde}
\left\{\begin{array}{l}
\hbox{For }p\in\mathcal{D} \hbox{ find }u\in X(\bar\Omega)\hbox{ s.t.}\\
a(u,v;p)=f(v;p,\phi),\quad \forall v\in X(\bar\Omega),
\end{array}\right.
\end{equation}
%
%\myc{MH1: Question: Do we want to consider an abstract setting for which we impose all assumptions which
%in general are used for the abstract RB setting? Then we simply could say that our forms satisfy the
%suppositions for the parameter separability which we need for the Min-Theta approach. We here
%could use a reference to e.g. \cite{H17}  As example we state and use the electrical machine with
%the forms from [24], which we also should (correctly) recall in this manuscript.\\
%%
%MH3: Do we need Lipschitz continuity of the forms w.r.t. parameter? In general we for our abstract problem
%simply should assume all requirements for the parameter separability, see f.e. the work of Haasdonk.
%For our application problem we then check the suppositions (which are satisfied according to [24].}
where the assumptions in problem \eqref{pde_par} hold true. The bilinear form $a(\cdot,\cdot; p)$ can be expressed with an affine linear decomposition:
\begin{equation}\label{affine}
a(u,v;p)=\sum_{q=1}^{Q_a} \sum_{i,j=1}^2  \hat\Phi_{q}^{i,j}(p)a_q^{i,j}(u,v),
\end{equation}
such that $\hat\Phi_{q}^{i,j}:\mathcal{D}\rightarrow\R\hbox{ for }q=1,\ldots, Q_a, i,j=1,2,$ is a function depending on $p$ and $a_q^{i,j}$ are the parameter independent continuous bilinear forms. In our example we have
$$
a_q^{i,j}(w,v):=\int_{\Omega_q(\bar p)}  \dfrac{\partial w}{\partial x_i} \dfrac{\partial v}{\partial x_j} \,\mathrm dx, \quad\mbox{for}\quad q = 1,\ldots, \blue{Q_a} \mbox{ and } i,j=1,2.
$$
%where $\bar p$ denotes some reference parameter.
%
%where $p^*$ is the reference parameter associate to the reference geometry $\bar{\Omega}$ and \myc{MH: i.e. $a_q^k$ is degenerate? With this definition we will not be able to guarantee uniform coercivity.}
%%
%$$
%\blue{
%\bar{\nu}^1(\blue{\phi}) = \left(\begin{array}{cc}
%                                  \nu^1(\phi) & 0\\
%                                  0 & 0
%                                 \end{array}
%\right)
%\quad\mbox{and}\quad
%\bar{\nu}^2(\blue{\phi}) = \left(\begin{array}{cc}
%                                  0 & 0 \\
%                                  0 & \nu^2(\phi)
%                                  \end{array}
%\right).
%}
%$$
%
\blue{This allows us to improve the computational efficiency in the evaluation of $a(u,v;p)$; the discrete approximation of the components $a_q^{i,j}(u,v)$ can be then computed once.} % and then are stored in the so called offline stage of the method. 
The same affine decomposition is applied to $f$ with
\begin{equation}\label{source_aff}
f(v;p,\phi)=\sum_{q=1}^{Q_f}  h_q(\phi)\hat\Phi_q^f(p) f_q(v)
\end{equation}
\blue{with $h_q(\phi)$ being} nonlinear differentiable functions depending on the uncertain parameter and
$$
f_q(v):=\int_{\Omega_q(\bar p)} \blue{(}f\blue{\circ\mathcal{T}_q}\blue{)} v \blue{|\det J_{\mathcal{T}_q}|} \,\mathrm dx, \quad\mbox{for}\quad q = 1,\ldots,Q_f.
$$
The assumption that the problem is affine dependent on the parameter $p$ is the key for the purpose of model reduction as we will see in Section \ref{sec:pmor}. We note that without loss of generality, in the present work, we consider the same number of subdomains $Q=Q_a=Q_f$ for the (bi)linear forms $a$ and $f$.

%\blue{We note that due the
%assumption of isotropic material we only have two components for $a_q^k$. Finally, we also note that without loss of generality, in the present work, we consider the same number of subdomains $Q$ for $a$ and $f$.}
Furthermore, we define the coercivity constant that will be a key ingredient in the certification of the model reduction algorithm, by
\begin{equation}\label{coer}
\alpha(p):=\inf_{w\in X(\bar\Omega), w\neq 0} \dfrac{a(w,w;p)}{\|w\|_{X(\bar{\Omega})}^2},
\end{equation}
and the continuity constant as
\begin{equation}\label{cont:cost}
\gamma(p):=\sup_{v\in X(\bar\Omega), v\neq 0}\sup_{w\in X(\bar\Omega), w\neq 0} \dfrac{a(w,v;p)}{\|w\|_{X(\bar{\Omega})}\|v\|_{X(\bar{\Omega})}}.
\end{equation}
To ease the notations, from here onwards, we drop the dependence on $\phi$ (e.g. $h_q(\phi) := 1$ in \eqref{source_aff}).

For the purpose of the optimization problem we will compute the
sensitivities ${u}^1_{i}: = \frac{\partial u(p)}{\partial p_i}$ which are
obtained by the derivative with respect to the parameters from equation
\eqref{pde} leading the following linear problem:
\begin{equation}\label{rid:sens}
\left\{\begin{array}{l}
\hbox{For }p\in\mathcal{D}\hbox{ find }u^1_{i}\in X(\bar{\Omega)}\hbox{ s.t.}\\
a(u^1_{i},v; p)=\frac{\partial f}{\partial p_i}(v;p)-\frac{\partial a}{\partial p_i}(u,v;p),\quad \forall v\in X(\bar\Omega),\, i=1,\ldots, N_p.
\end{array}\right.
\end{equation}
We note that due to the affine decomposition the $i$-th partial derivatives
of $a(\cdot,\cdot;p)$ and $f(\cdot;p)$ are given by the derivatives of
coefficient functions $\hat\Phi_{q}^{i,j}$ and $\hat\Phi_q^f$, $q = 1,\ldots,Q$, $i,j=1,2$ and
can be computed analytically. More generally, for a multiindex
$\mu\in\N_0^{N_p}$, $|\mu|=n$,
the $\mu$-th sensitivity can be
computed by the following proposition.

\begin{proposition}
Let the coefficient functions $\hat\Phi_{q}^{i,j}(p)$ and $\hat\Phi_q^f(p,\phi)$, $q =
1\ldots,Q$, $i,j=1,2$ be $n$-times differentiable with respect to $p$.  Then the
solution $u$ is differentiable with respect to $p$ and the sensitivities
$u_\mu^n = \frac{\partial^\mu u}{\partial p^\mu} \in \blue{X(\bar\Omega)}$,
$|\mu|=n$, satisfy the sensitivity equation

 \begin{equation}
 \label{eq:sens_general}
 a(u_\mu^n,v;p) = \frac{\partial^\mu f}{\partial p^\mu}(v;p) 
- \sum_{1\le |\kappa|\le n} \binom{\mu}{\kappa}
\frac{\partial^\kappa a}{\partial p^\kappa}(u^{n-|\kappa|}_{\mu-\kappa},v;p)\quad 
\forall v \in \blue{X(\bar\Omega)},
\end{equation}
where $\binom{\nu}{\kappa} = \prod_{i=1}^{N_p}\binom{\mu_i}{\kappa_i}=
\frac{\mu!}{\kappa!(\mu-\kappa)!}$ denotes the binomial coefficient
for multiindices.
\end{proposition}

The proof of the proposition follows from the general Leibniz rule for the
$\mu$-th derivative applied to \eqref{affine} and is omitted here.
We can easily see that for
$n=0$ we find the state equation \eqref{pde} and for $n=1$ the first
sensitivity equation \eqref{rid:sens}.

% \begin{rmk}
%  In the case the parameter $p$ enters \eqref{pde} only through the linear form $f(v;p)$, i.e., the bilinear form is independent of $p$, the sensitivity equation reduce to 
%  %
%  \begin{equation}
%  \label{eq:sens_general_rhs}
%  a(u^n,v) = \frac{\partial^n f}{\partial p^n}(v;p).
%  \end{equation}
%  %
% \end{rmk}

Next, we introduce an high dimensional finite element (FE) discretization of our model problem in the space $X_h\blue{(\bar\Omega)}\subset X\blue{(\bar\Omega)}$. The discrete problem then reads
\begin{equation}\label{rid:fem}
\left\{\begin{array}{l}
\hbox{For }p\in\mathcal{D}\hbox{ find }u_h\in X_h(\bar{\Omega})\hbox{ s.t.}\\
a(u_h,v_h; p)=f(v_h;p),\quad \forall v_h\in X_h(\bar\Omega).
\end{array}\right.
\end{equation}
For the discrete problem we use the ansatz $u_h = \sum_{i =1}^{N} ({\bf u}_h)_i \varphi_i$, where $\varphi_i$ are suitable FE ansatz functions. Problem \eqref{rid:fem} then is equivalent to the linear system
\begin{equation}\label{rid:fem2}
\left\{\begin{array}{l}
\hbox{For }p\in\mathcal{D}\hbox{ find } {\bf u}_h(p)\in \mathbb R^N\hbox{ s.t.}\\
{\bf K}(p){\bf u_h}(p)={\bf f}(p),
\end{array}\right.
\end{equation}
where ${\bf K}(p)\in\R^{N\times N}$ is the stiffness matrix $\left({\bf K}(p)\right)_{ij}=a(\varphi_j,\varphi_i;p)$, $1\le i,j \le N$, and the right hand side ${\bf f}(p)\in\R^N$ is obtained by $\left({\bf f}(p)\right)_i=f(\varphi_i;p)$, $1 \le i \le N$. The FE system matrix keeps the dependency on the parameter $p$ and we apply the affine decomposition to it in the following way 
\begin{equation}
\label{eq:affine}
{\bf K}(p) = \sum_{q = 1}^{Q}\sum_{i,j=1}^2 \hat\Phi_q^{i,j}(p) {\bf K}_q^{i,j},
\end{equation}
where ${\bf K}_q^{i,j}$,  $q = 1,\ldots,Q$ and $i,j=1,2$, are the system matrices on the $Q$ sub-domains.
% todo Q is not no of subdomains, since p=(k,i,j), k subdomain index
We note again that the $p$ dependency of ${\bf K}$ appears only in the weight functions $\hat\Phi_q^{i,j}$ which are easy and inexpensive to evaluate.
The same applies for the right hand side. Similarly, the discretized version of equation \eqref{rid:sens} reads
\begin{equation}
\label{eq:sens}
 {\bf K}(p){\bf u}^1_{h,i}(p) = {\bf \tilde{f}}^1_i, \quad\mbox{for }i = 1,\ldots,N_p,
\end{equation}
with
$$
{\bf \tilde{f}}^1_i = \frac{\partial{\bf f}}{\partial p_i}(p) - \frac{\partial{\bf K}}{\partial p_i}(p){\bf u}_h \quad\mbox{for }i = 1,\ldots,N_p,
$$
%
 %$N_p$ is the number of parameters in the problem and
 \blue{where the subindex $i$ indicates as above the derivative with respect to the $i$-th parameter
$p_i$.}  Note
that these derivatives are easy to compute due to the previously introduced
affine decomposition \eqref{eq:affine}.  The derivative of the matrix ${\bf
K}(p)$ and the vector ${\bf f}(p)$ are given by the derivatives of the
functions $\hat\Phi_q^{i,j}(p)$ and $\hat\Phi_q^f(p,\phi)$, $q = 1,\ldots,Q$, $i,j=1,2$,
respectively. For the general sensitivity equations \eqref{eq:sens_general}
we get with a multiindex $\mu\in\N_0^{N_p}$, $|\mu|=n$,
\begin{equation}\label{Neq:sens}
{\bf K}(p){\bf u}^n_{h,\mu}=\dfrac{\partial^\mu {\bf f}}{\partial p^\mu}(p)
-\sum_{1\le |\kappa|\le n} \binom{\mu}{\kappa}\dfrac{\partial^\kappa {\bf K}}{\partial
p^\kappa}(p){\bf u}^{n-|\kappa|}_{h,\mu-\kappa}.
\end{equation}
In the rest of the work we will continue to use the following compact notation when appropriate:
$${\bf u}_{h,\mu}^{n}(p):=
\dfrac{\partial^\mu {\bf u}_h}{\partial p^\mu}(p)$$
in order to denote the $\mu$-th derivative of ${\bf u}_h$ with respect
to $p$ for a multindex $\mu\in \N_0^{N_p}$, $|\mu|=n$.
Moreover, we sometimes use the fact that the state is the zeroth order
sensitivity, i.e., ${\bf u}_{h,\mu}^{0}(p)={\bf u}_{h}(p)$ with
$|\mu|=0$.
Similarly, we will adopt ${\bf f}^n_{\mu}$ and ${\bf K}^n_{\mu}$.

From now on we will focus on the discrete version (\ref{rid:fem2}) of the problem (\ref{pde_par}). All further steps are analogous in the continuous setting.

% -------------------------------------------------
% -------------------------------------------------
% --------------- SECTION 3 -----------------
% -------------------------------------------------
% -------------------------------------------------
\section{Optimization problem} 
\label{sec:opt}
This section is devoted to the optimization and the robust optimization problem. We will formulate an optimization problem governed by a parametrized PDE, then investigate the nominal optimization and its robust counterpart. By the robust formulation, uncertainties in model parameters are taken into account. This is done by utilizing a worst-case formulation. For an efficient realization, approximation techniques of different orders are investigated. Here we closely follow \cite{LU17} and repeat the ingredients for the convenience of the reader.

\subsection{The nominal optimization}
\label{sec:opt_nominal}
We first introduce the nominal optimization problem. In this setting all model parameters $\phi$ are fixed to a nominal value $\hat\phi$. The optimization problem then reads as:
\begin{equation}
 \label{eq:opt}
 \begin{array}{rl}
   \displaystyle\min_{p\in \mathbb R^{N_p}, {\bf u}_h \in \mathbb R^{N}} & \tilde g_0({\bf u}_h,p,\hat \phi),\\
   \mbox{subject to (s.t.)} & \tilde g_i({\bf u}_h,p,\hat \phi) \le 0,\quad i = 1,\ldots,N_g,\\
                            & e({\bf u}_h,p,\hat \phi) = 0,
 \end{array}
\end{equation}
where $p\in\mathbb R^{N_p}$ is the optimization variable and $\hat
\phi \in \mathcal U_k\subset \mathbb R^{N_\phi}$ is the fixed model
parameter.  The function $\tilde g_0$ is the objective function and $\tilde
g_i$ are the $N_g\in \mathbb R$ inequality constraints with $\tilde g_i:
\mathbb R^{N}\times\mathbb R^{N_p}\times\mathbb R^{N_\phi} \mapsto \mathbb
R$ for $i = 0,\ldots,N_g$.  Furthermore, $e$ is an equality constraint
governed by the discretized parametric PDE in \eqref{rid:fem2}, i.e.,
$e(\blue{\bf u}_h,p,\hat \Phi) = 0 \Longleftrightarrow K(p,\hat \Phi)u_h=f(p,\hat
\Phi)$.  We assume $\tilde g_i$, $i = 0,\ldots,N_g$, and $e$ to be
continuously differentiable.  It directly follows from the coercivity of the
bilinear form $a$ that the parametric PDE for every admissible parameter $p$
admits a unique solution $u_h$, and that the Jacobian $\frac{\partial
e}{\partial \blue{\bf u}}$ is boundedly invertible, compare for example \cite{HPUU09}. 
Further, let $e({\bf u}_h,p,\phi) = 0$ have a unique solution ${\bf u}_h =
{\bf u}_h(p)$ for every admissible $p$.  Then we can introduce the reduced
objective and constraint functions

%\myc{MH4: Should we be more specific about the dependance on $\Phi?$ Maybe one could add a sentence on what $\Phi$ in the example is.}
%
$$
g_i(p,\hat\phi) := \tilde g_i({\bf u}_h,p,\hat\phi)\quad\mbox{for}\quad i= 0,\ldots,N_g.
$$
The reduced optimization problem associated with \eqref{eq:opt} then reads
\begin{equation}
  \label{eq:opt_red}
  \begin{array}{rl}
    \displaystyle\min_{p\in \mathbb R^{N_p}} & g_0(p,\hat\phi),\\
    \mbox{s.t.} & g_i(p,\hat\phi) \le 0,\quad i = 1,\ldots,N_g.
   \end{array}
\end{equation}
Note that the reduced objective and the reduced constraints $g_i(p)$, $i = 0,\ldots,N_g$, are continuously differentiable due to the previous assumptions. Since $e$ is continuously differentiable with invertible Jacobian $\blue{\frac{\partial}{\partial {\bf u}}e({\bf u},p,\phi)}$, the implicit function theorem guarantees that also ${\bf u}_h$  is continuously differentiable with respect to the parameter. In the further investigation, we focus on the reduced formulation.

\subsection{The robust optimization}
\label{sec:opt_robust}
After having introduced the nominal optimization problem we will now introduce the robust optimization framework, where the model parameter $\phi$ is considered to be uncertain. In the presented setting, the uncertain parameter is a model parameter but the framework is not limited to this setting. We will next formulate a robust version of \eqref{eq:opt_red}. For this, we first have to make some assumptions on the uncertain parameter and introduce the {\it uncertainty set}
\begin{equation}
\label{un_set}
 \mathcal U_k = \left\{\phi\in\mathbb R^{N_\phi}\,\big|\, \|D^{-1}(\phi - \hat\phi)\|_k \le 1\right\},
\end{equation}
where $\hat\phi$ is a nominal value and $D\in\mathbb R^{N_\phi\times N_\phi}$ an invertible scaling matrix. We assume that the uncertain parameter remains within the given bounded set $\mathcal U_k$. For different choices of $k$ we get different uncertainty sets. The two most commonly used ones are $k = \{2,\infty\}$. The case $k = \infty$ is a special case and allows a representation using upper and lower bounds \blue{$\overline{\phi}, \underline{\phi}$, respectively}. For this we set $D = \mathrm{diag}((\overline{\phi} - \underline{\phi})/2)$ and $\hat\phi = (\underline{\phi} + \overline{\phi})/2$. The uncertainty set can then be written as
$$
\mathcal U_\infty = \{\phi\in\mathbb R^{N_\phi}\,\big|\,\underline{\phi}\le\phi\le\overline{\phi}\}.
$$
Next we want to focus on the robust optimization problem utilizing the worst-case formulation \cite{BN02,DBK06,Zha07}. Using the uncertainty set we can define the {\it worst-case function} as
$$
 \varphi_i(p) := \max_{\phi\in \mathcal U_k} g_i(p,\phi),\quad i=0,\ldots,N_g.
$$
For an interpretation, for every fixed parameter $p$ the function $\varphi_i: \mathcal U_k \rightarrow \mathbb R$ is the worst-case value of the function $g_i$ with $\phi\in \mathcal U_k$. Using the definition of the uncertainty set, the worst-case function can be rewritten as
\begin{equation}
\label{eq:worstcasefunction}
  \begin{array}{rl}
\varphi_i(p) := \displaystyle\max_{\phi\in \mathbb R^{N_\phi}} & g_i(p,\phi)\\ 
\mbox{s.t.}\,\, & \|D^{-1}(\phi - \hat\phi)\|_k \le 1,\quad i=0,\ldots,N_g.
  \end{array}
\end{equation}
Using the worst case function we will formulate the robust optimization problem. The {\it robust counterpart} of \eqref{eq:opt_red} is given by
\begin{equation}
 \label{eq:robust}
 \begin{array}{rl}
  \displaystyle\min_{p\in\mathbb R^{N_p}} & \varphi_0(p)\\
  \mbox{s.t.}\,\, & \varphi_i(p) \le 0,\quad i=1,\ldots,N_g.
 \end{array}
\end{equation}
A solution $p$ of \eqref{eq:robust} is robust against uncertainties in the parameter $\phi$ and is hence referred to as {\it robust optimal solution}. Note that the solution is feasible for \eqref{eq:opt_red} for all $\phi\in\mathcal U_k$ and optimal with respect to the objective function.

The introduced robust optimization problem \eqref{eq:robust} is of bi-level structure and hence difficult to solve. Thus it is required to develop tailored methods to solve the problem. For general nonlinear problems in \cite{BN02}, it is proposed to employ approximations. The {\it approximated robust counterpart} of \eqref{eq:robust} is then introduced as
\begin{equation}
 \label{eq:robust_approx}
 \begin{array}{rl}
  \displaystyle\min_{p\in\mathbb R^{N_p}} & \hat \varphi_0(p)\\
  \mbox{s.t.}\,\,&\hat\varphi_i(p) \le 0,\quad i=1,\ldots,N_g,
 \end{array}
\end{equation}
where $\hat\varphi_i$ can be computed more efficiently compared to \eqref{eq:worstcasefunction} and is referred to as the {\it approximated worst-case function}. First and second order approximations $\hat\varphi_i$ of the worst-case function $\varphi_i$ will be considered. The first order approximation for the general nonlinear case was investigated in \cite{DBK06,Zha07} while the second order approximation is a recent approach \cite{LU17} and can be seen as a modification of \cite{Sic13}.

\subsubsection{Linear approximation of the robust counterpart}
\label{sec:opt_robust_linear}
In the first order approach suggested in \cite{DBK06,Zha07}, a linearization of the worst-case function is carried out. For this a nominal value $\hat \phi$ for the uncertain parameter is chosen and the approximated worst-case function is then given by the first order Taylor expansion
$$
\begin{array}{rl}
\hat\varphi_i(p) := \displaystyle\max_{\delta_i\in\mathbb R^{N_\phi}} & g_i(p,\hat\phi) + \nabla_\phi g_i(p,\hat\phi)^\top\delta_i\\
\mbox{s.t.}\,\,& \|D^{-1}\delta_i\|_k\le 1,\quad i=0,\ldots,N_g,
 
\end{array}
$$
where $\delta_i = \phi - \hat\phi$. The solution of the resulting convex optimization problem can be given analytically in the form
\begin{equation*}
 \label{eq:worstcasefunction_linear}
 \hat\varphi_i(p) = g_i(p,\hat\phi) + \|D \nabla_\phi g_i(p,\hat\phi)\|_{k^*}
\end{equation*}
with $\|D \cdot\|_{k^*}$ the dual norm of $\|D^{-1}\cdot\|_{k}$ for $k^* = k/(k-1)$ with $k^* = 1$ for $k = \infty$. Note that this is a standard result for scaled H\"older norms. Utilizing these results, the linear approximated robust counterpart reads as
\begin{equation}
\label{eq:robust_approx_linear}
\begin{aligned}
 &\min_{p\in\mathbb R^{N_p}} g_0(p,\hat\phi) + \|D\nabla_\phi g_0(p,\hat\phi)\|_{k^*}\\
 &\quad\mbox{s.t.}\,\, g_i(p,\hat\phi) + \|D\nabla_\phi g_i(p,\hat\phi)\|_{k^*} \le 0,\quad i=1,\ldots,N_g.
\end{aligned}
\end{equation}
We refer to \eqref{eq:robust_approx_linear} as the {\it linear approximated robust counterpart} of \eqref{eq:opt_red}. Due to the norm the objective function and the inequality constraints are non-differentiable if the term inside the norm becomes zero. For the case $k = \infty$, the remedy is to introduce slack variables, i.e.,
$$
\begin{array}{rlr}
 \displaystyle\min_{p\in\mathbb R^{N_p},\zeta_i\in\mathbb R} & g_0(p,\hat\phi) + \zeta_0 &\\
 \mbox{s.t.}\,\, &g_i(p,\hat\phi) + \zeta_i \le 0,&i=1,\ldots,N_g,\\
 &-\zeta_i \le D\nabla_\phi g_i(p,\hat\phi) \le \zeta_i,& i=0,\ldots,N_g.
\end{array}
$$
For the case $k = 2$, the problem \eqref{eq:robust_approx_linear} can
be reformulated by using second order cone constraints that can be handled
effciently for example by interior point techniques.
But usually the solution of \eqref{eq:robust_approx_linear} is
sufficiently far away from points of nonsmoothness and in practice
standard algorithms for smooth nonlinear problems work well
without any modification. Nevertheless, one has to verify that the robust optimal solution does not lie in a non-differentiable point \cite{DBK06,Sic13}.

The required derivatives for the approximations in this approach and the optimization can be computed using either the adjoint or the sensitivity approach \cite{HPUU09}. Note that the adjoint approach is desirable when the number of uncertain parameters is large. For a detailed discussion about the different approaches we refer the reader to \cite{DBK06}.

While computationally attractive, the first order approximation can suffer from an inaccurate approximation. Hence the influence of the uncertain parameter $\phi$ might be described insufficiently, which was already observed in \cite{DGMM08}. This leads to the development of second order methods \cite{LU17,Sic13}.

\subsubsection{Quadratic approximation of the robust counterpart}
\label{sec:opt_robust_quadratic}
In the second order approach, a quadratic approximation of the worst case function is utilized. This idea has been investigated in \cite{LU17,Sic13}. We will apply the method directly to the reduced problem \eqref{eq:opt_red}. Since this approach is recent, we will give a short overview of the strategy.

The second order approximation is generated using the second order Taylor expansion of \new{$g_i(p,\cdot)$ around $\hat\phi$, i.e.,}
\begin{equation}\new{
    g_i(p,\hat\phi + \delta_i)
    \approx
    q_i(p,\hat\phi,\delta_i)
    \coloneqq
    \alpha_i(p,\hat\phi)
    +
    c_i(p,\hat\phi)^\top\delta_i
    +
    \tfrac{1}{2}
    \delta_i^\top H_i(p,\hat\phi)\delta_i,}
\end{equation}
\new{where $\alpha_i(p,\hat\phi) = g_i(p,\hat\phi)$, $c_i(p,\hat\phi) = \nabla_\phi g_i(p,\hat\phi)$, and $H_i(p,\hat\phi) = \nabla_{\phi\phi} g_i(p,\hat\phi)$.
As before, we have set $\delta_i = \phi - \hat\phi$.
The corresponding approximated worst-case function is}
\begin{equation}
 \label{eq:trustregion}
 \hat\varphi_i(p) := \new{\max_{\delta_i\in\mathbb R^{N_\phi}} q_i(p,\hat\phi,\delta_i)}
 \quad\mbox{s.t.}\quad \|D^{-1}\delta_i\|_k\le 1,\, i=0,\ldots,N_g.
\end{equation}
For the quadratic approximation we do not have a closed solution as in the linear case. Fortunately, the problem exhibits a well known structure. In the case $k = 2$, the problem corresponds to a trust region problem and
is well-studied. The exact solutions \blue{of} \eqref{eq:trustregion} are characterized
as follows.

\begin{theorem}
\label{thm:trs}
 \new{For given $p$, the} vector $\delta_i\in\mathbb R^{N_\phi}$ is a global solution of the trust region problem \eqref{eq:trustregion} if and only if there exists a Lagrange multiplier $\lambda_i\in\mathbb R$ satisfying
 $$
 \begin{array}{c}
  \lambda_i \ge 0,\quad 
  \|D^{-1}\delta_i\|_2 \le 1,\quad
  \lambda_i(\|D^{-1}\delta_i\|_2 -1) = 0,\\[1em]
  (-H_i(\new{p},\hat\phi) + \lambda_i \mathbb D)\delta_i = c_i(\new{p},\hat\phi)
 \end{array}
 $$
with $-H(\new{p},\hat\phi) + \lambda_i\mathbb D$ positive semidefinite and $\mathbb D = D^{-\top}D^{-1}$.
\end{theorem}

The proof of this theorem can be found in \cite{CGT00,HPUU09}. By adding a square to the norms in the constraints, we obtain the \new{equivalent but} differentiable formulation of the {\it quadratic approximated robust counterpart} of \eqref{eq:opt_red} by
%
% \begin{equation}
%  \label{eq:robust_approx_quadratic}
%  \begin{aligned}
%  \min_{\substack{
%         p\in\mathbb R^{N_p},\\
%         \delta_0,\ldots,\delta_{N_g}\in\mathbb R^{N_\phi},\\
%         \lambda_0,\ldots,\lambda_{N_g}\in\mathbb R}}
%         g_0(p,\hat\phi) + \nabla_\phi g_0(p,\hat\phi)^\top\delta_0 + \frac{1}{2}& \delta_0^\top H_0(p,\hat\phi)\delta_0\\
%  \quad\mbox{s.t.}\quad g_i(p,\hat\phi) + \nabla_\phi g_i(p,\hat\phi)^\top\delta_i + \frac{1}{2} \delta_i^\top H_i(p,\hat\phi)\delta_i \le 0,&\quad i = 1,\ldots,N_g,\\
%  \left(\begin{array}{c}
%         -\nabla_{\phi}g_i(\cdot,\hat\phi)-H_i(\cdot,\hat\phi)\delta_i + \lambda_i \mathbb D \delta_i\\
%         \lambda_i(\|D^{-1}\delta_i\|_2^2 -1)
%        \end{array}\right) = 0,&\quad i = 0,\ldots,N_g,\\
%  \|D^{-1}\delta_i\|_2^2 -1 \le 0,&\quad i = 0,\ldots,N_g,\\
%  -\lambda_i \le 0,&\quad i = 0,\ldots,N_g,\\
%  H_i(\cdot,\hat\phi) - \lambda_i \mathbb D \preceq 0,&\quad i = 0,\ldots,N_g,
%  \end{aligned}
% \end{equation}
% uaie
\begin{equation}
 \label{eq:robust_approx_quadratic}
 \begin{aligned}
 \min\quad & \new{\alpha_0(p,\hat\phi) + c_0(p,\hat\phi)^\top\delta_0 + \tfrac{1}{2} \delta_0^\top H_0(p,\hat\phi)\delta_0}\\
 \text{s.t.}\quad & \begin{alignedat}[t]{2}
    \new{\alpha_i(p,\hat\phi) + c_i(p,\hat\phi)^\top\delta_i + \tfrac{1}{2} \delta_i^\top H_i(p,\hat\phi)\delta_i} &\le 0,
    &\quad& i = 1,\ldots,N_g,\\
 \left(\begin{array}{c}
        \new{-c_i(p,\hat\phi)}-H_i(\new{p},\hat\phi)\delta_i + \lambda_i \mathbb D \delta_i\\
        \lambda_i(\|D^{-1}\delta_i\|_2^2 -1)
       \end{array}\right) &= 0,&\quad& i = 0,\ldots,N_g,\\
 \|D^{-1}\delta_i\|_2^2 -1 &\le 0,&\quad& i = 0,\ldots,N_g,\\
 -\lambda_i &\le 0,&\quad& i = 0,\ldots,N_g,\\
 H_i(\new{p},\hat\phi) - \lambda_i \mathbb D &\preceq 0,&\quad& i = 0,\ldots,N_g,
    \end{alignedat}
 \end{aligned}
\end{equation}
where $A \preceq 0$ denotes that $A$ is a negative semidefinite matrix.
\new{The optimization variables are $p\in\mathbb R^{N_p}$, $\delta_0,\ldots,\delta_{N_g}\in\mathbb R^{N_\phi}$, and $\lambda_0,\ldots,\lambda_{N_g}\in\mathbb R$.}
The semidefinite constraint can be reformulated using the smallest
eigenvalues as was outlined in \cite{Sic13} or be treated, e.g., by interior
point techniques.  Problem
\eqref{eq:robust_approx_quadratic} is a {\it mathematical program with
equilibrium constraints} (MPEC).  To solve these problems numerically one
has to pay attention to the complementarity constraint.  It turns out that
SQP method under relatively mild assumptions and few modifications
\cite{FLRS06,Ley06} can solve these type of problems efficiently.  In
\cite{LU17} the required changes are outlined in details for the given
problem.

The computation of the derivatives required for the quadratic
approximation can be done by different approaches.  Depending on the number
of parameters, the sensitivity or the adjoint approach should be chosen.  In
this work, it is assumed that the number of uncertain parameters is small
and the sensitivity based approach is utilized.

\subsubsection{\new{Improving the approximation by moving the expansion point}}
\label{ssec:shifting}

The quadratic approximation is often a notable improvement over the linear approximation.
Nonetheless, Taylor approximations generally are accurate only locally, that is, in a neighborhood of the expansion point $\hat\phi$, and so the second-order model can sometimes be a poor predictor for the effects of the uncertain parameters, especially when the uncertainty set is relatively large and the expansion point is far away from the worst case.
In these cases, the accuracy can be improved by moving the expansion point of the Taylor approximation towards a presumed maximum of the modeled function. 
This can be implemented with the aid of an additional set of variables $\hat\phi_i \in \R^{N_\phi}$, $i = 0,\ldots,N_g$, for the individual expansion points of the models.
The second-order Taylor expansion of $g_i(p,\dummy)$ around $\hat\phi_i$ evaluated near the center $\hat\phi$ of the uncertainty set is given by
\begin{align*}
    g_i(p,\hat\phi + \delta_i)
    \approx
    q_i^s(p,\hat\phi_i,\delta_i)
    \coloneqq{}& g_i(p,\hat\phi_i)
     + \nabla_\phi g_i(p,\hat\phi_i)^\top(\delta_i - \hat\delta_i)
    \\
    &\quad
     + \tfrac{1}{2}(\delta_i - \hat\delta_i)^\top \nabla_{\phi\phi}g_i(p,\hat\phi_i)(\delta_i - \hat\delta_i)
    \\
    ={}& \alpha_i^s(p,\hat\phi_i) + c_i^s(p,\hat\phi_i)^\top\delta_i + \tfrac{1}{2}\delta_i^\top H_i(p,\hat\phi_i)\delta_i,
\end{align*}
where $\hat\delta_i = \hat\phi_i - \hat\phi$ as well as $H_i(p,\hat\phi_i) = \nabla_{\phi\phi} g_i(p,\hat\phi_i)$, and additionally
\begin{align*}
    \alpha_i^s(p,\hat\phi_i) &= g_i(p,\hat\phi_i) - \nabla_\phi g_i(p,\hat\phi_i)^\top\hat\delta_i,
    \\
    c_i^s(p,\hat\phi_i) &= \nabla_\phi g_i(p,\hat\phi_i) - \nabla_{\phi\phi} g_i(p,\hat\phi_i)\hat\delta_i.
\end{align*}
Since a good choice for $\hat\phi_i$ is often not known a priori, it can be adaptively chosen in the course of the optimization.
A schematic algorithm that does so is stated in Algorithm~\ref{Alg:SQP}, which is based on a basic version of the SQP method by Powell, see~\cite{P83}.
\begin{algorithm}
\caption{(Schematic adaptive SQP algorithm)}
\label{Alg:SQP}
\begin{algorithmic}[1]
\REQUIRE Initial guess $(p^0,(\delta_i^0)_i,(\lambda_i^0)_i)$ and $(\hat\phi_i^0)_i$, penalty parameter $\rho > 0$
\WHILE{termination criterion not met}
\STATE Solve a quadratic program with expansion points $(\hat\phi_i^k)_i$ to obtain update $(\Delta p^k,(\Delta \delta_i^k)_i,(\Delta \lambda_i^k)_i$ \label{step:qp}
\STATE Determine step length $\alpha_k$ to achieve proper descent in $P_\rho^1$ \label{step:descent}
\STATE Set $(p^{k+1},(\delta_i^{k+1})_i,(\lambda_i^{k+1})_i) = (p^{k},(\delta_i^{k})_i,(\lambda_i^{k})_i) + \alpha_k\cdot(\Delta p^k,(\Delta \delta_i^k)_i,(\Delta \lambda_i^k)_i)$
\STATE Determine update $\Delta\hat\phi_i^k$ to increase accuracy of $q_i^s(p^{k+1},\hat\phi_i^k,\dummy)$, $i = 0,\ldots,N_g$
\STATE Set $\hat\phi_i^{k+1} = \hat\phi_i^k + \Delta\hat\phi_i^k$, $i = 0,\ldots,N_g$
\STATE $k \gets k + 1$
\ENDWHILE
%\STATE Compute ${\bf u}^{i,\ell}(p^{k+1})$ and ${\bf u}^i(p^{k+1})$, 
%\ENDWHILE
\end{algorithmic}
\end{algorithm}
%
%Since globally maximizing $g_i(p,\dummy)$ whenever $p$ is updated would be too expensive, it is reasonable to move the expansion points only gradually.
The algorithm generates a sequence of optimization variables $x^k \coloneqq (p^{k},(\delta_i^{k})_i,(\lambda_i^{k})_i)$ and a sequence of expansion points $\hat\Phi^k \coloneqq ((\hat\phi_i^k)_i)$, where $k \in \N$ is the iteration counter.
The feasibility is eventually enforced with an exact $\ell_1$-penalty function, denoted here by $P_\rho^1(x;\hat\Phi)$, with a penalty parameter $\rho > 0$.
Recall that, if we write~\eqref{eq:robust_approx_quadratic} in the form $\min_x\,\{\, f(x;\hat\Phi) : h_\leq(x;\hat\Phi) \leq 0,\, h_=(x;\hat\Phi) = 0 \,\}$ for brevity, the penalty function is defined by
\[
    P_\rho^1(x;\hat\Phi)
    =
    f(x;\hat\Phi)
    +
    \rho\cdot\lVert \max\{0, h_\leq(x;\hat\Phi) \} \rVert_1
    +
    \rho\cdot\lVert h_=(x;\hat\Phi) \rVert_1,
\]
where the maximum is to be understood componentwise, and $\lVert\dummy\rVert_1$ denotes the $1$-norm.
We assume that a sufficiently large $\rho$ is known in advance, but in practise it is usually chosen dynamically.

The quadratic program in line~2 is obtained from~\eqref{eq:robust_approx_quadratic} at the current iterate in the standard way, i.e., by linearization of the constraints and by a quadratic model of the objective function with a positive definite approximation of the Hessian of the associated Lagrangian.
The novelty is that the problem functions of~\eqref{eq:robust_approx_quadratic} are changed between iterations by replacing $\alpha_i(p,\hat\phi)$, $c_i(p,\hat\phi)$ and $H_i(p,\hat\phi)$ by $\alpha_i^s(p,\hat\phi_i^k)$, $c_i^s(p,\hat\phi_i^k)$ and $H_i(p,\hat\phi_i^k)$ before line~2.

There is some freedom in how the updates $\Delta\hat\phi_i^k$ are determined.
One option is to perform a single step of a projected gradient method to maximize $g_i(p^k,\hat\phi_i^k + \dummy)$.
In this case, we choose $\hat\phi_i^{k+1}$ to be the projection of $\hat\phi_i^k + \sigma_k\nabla_\phi g_i(p^k,\hat\phi_i^k)$ onto the uncertainty set, where $\sigma_k > 0$ is an adequate step length.
For details, we refer to Section~2.2.2 of~\cite{HPUU09}.
This strategy tries to increase the accuracy of the quadratic models by gradually moving the expansion points towards a local maximum of the modeled functions.

Since we change the MPEC formulation~\eqref{eq:robust_approx_quadratic} in between two SQP iterations, standard convergence results do not apply.
The major issue is that we cannot use the common monotonicity argument about the sequence of penalty function values, which is used to show that the penalty function values converge.
A simple way to fix this is given in the next lemma.

\begin{lemma}
\label{lemma:sqp}
    Assume that the sequence of expansion points $(\hat\Phi^k)_k$ converges as $k \to \infty$.
    Furthermore, assume that
    \[
        \Big\lvert \sum_{k=0}^\infty P_\rho^1(x^{k+1};\hat\Phi^{k+1}) - P_\rho^1(x^{k+1};\hat\Phi^k) \Big\rvert < \infty.
    \]
    If the sequence $(x^k)_k$ has an accumulation point, then
    \[
        \lim_{k\to\infty} P_\rho^1(x^{k+1};\hat\Phi^{k+1}) - P_\rho^1(x^{k};\hat\Phi^{k}) = 0.
    \]
\end{lemma}
\begin{proof}
    Let $x^*$ be an accumulation point of $(x^k)_k$ and let $K \subseteq \N$ be an infinite subset such that $x^k \xrightarrow{K} x^*$.
    The continuity of $P_\rho^1$ implies that $(P_\rho^1(x^k;\hat\Phi^k))_k$ is a Cauchy sequence.
    Hence, for any $\varepsilon > 0$ there is some $N(\varepsilon) \in \N$ such that
    \[
        \lvert P_\rho^1(x^{k_2};\hat\Phi^{k_2}) - P_\rho^1(x^{k_1};\hat\Phi^{k_1}) \rvert
        < \varepsilon,
        \quad
        \text{for all}
        \quad
        k_1,k_2 \in K,\; k_2 \geq k_1 \geq N(\varepsilon). 
    \]
    Noting $P_\rho^1(x^{k_2};\hat\Phi^{k_2}) - P_\rho^1(x^{k_1};\hat\Phi^{k_1}) = \sum_{k=k_1}^{k_2-1} P_\rho^1(x^{k+1};\hat\Phi^{k+1}) - P_\rho^1(x^{k};\hat\Phi^k)$, we obtain
    \begin{align*}
        \Big\lvert \sum_{k=k_1}^{k_2-1}
          P_\rho^1(x^{k+1};\hat\Phi^{k+1}) - P_\rho^1(x^{k+1};\hat\Phi^{k})
        + P_\rho^1(x^{k+1};\hat\Phi^k) - P_\rho^1(x^{k};\hat\Phi^k)
        \Big\rvert < \varepsilon.
    \end{align*}
    Since $P_\rho^1(x^{k+1};\hat\Phi^k) - P_\rho^1(x^{k};\hat\Phi^k) \leq 0$ by line~3 of Algorithm~\ref{Alg:SQP}, the summability assumption shows that $0 > \sum_{k=k_1}^{k_2-1} P_\rho^1(x^{k+1};\hat\Phi^k) - P_\rho^1(x^{k};\hat\Phi^k) > -\infty$ as $k_2 \xrightarrow{K} \infty$.
    This implies $P_\rho^1(x^{k+1};\hat\Phi^k) - P_\rho^1(x^{k};\hat\Phi^k) \to 0$ for $k \to \infty$, from which the claim follows.
\end{proof}

The lemma can be used to apply standard convergence results to Algorithm~\ref{Alg:SQP}.
For example, we can argue analogously to the proof of Theorem~1 in~\cite{P83} in order to show that every accumulation point of $(x^k)_k$ is a KKT point of~\eqref{eq:robust_approx_quadratic}, under the additional, usual assumptions of SQP methods like uniform boundedness of the Hessian approximations of the Lagrangian.
In this proof, the existence of a limit of the sequence $(P_\rho^1(x^{k};\hat\Phi^{k}))_k$, which is established by Lemma~\ref{lemma:sqp},
implies that every accumulation point of $(x^{k})_k$ is feasible.
The rest of the proof is simple to modify, so that we skip the details.

% -------------------------------------------------
% -------------------------------------------------
% --------------- SECTION 4 -----------------
% -------------------------------------------------
% -------------------------------------------------

\section{Proper Orthogonal Decomposition for Parametrized problems}\label{sec:pmor}
Numerical approximation for robust optimization problems can be expensive since it involves the solution of several PDEs. Furthermore, the sensitivity approach enlarges the number of PDEs and it increases the computational costs of its approximation. For this reason, in this section, we introduce a model order reduction technique to reduce the complexity of the problem.

Here, we focus on the POD method for the approximate solution of the
parametrized equation (\ref{rid:fem2}).  \blue{As before,} ${\bf
u}_h(p)\in\mathbb R^N$ is the model vector associated to the FE solution of
\eqref{rid:fem2} for a given parameter $p\in \mathcal{D}\subset \R^{N_p}$
 and ${\bf u}^n_{h,\mu}(p)\in\mathbb R^N$ \red{are the sensitivities} of order $n$
according to \eqref{Neq:sens}.
For this purpose let $\{p^j\}_{j=1}^m$ be some points in
$\mathcal{D}$ and let ${\bf u}_h(p^j), {\bf u}^n_{h,\mu}(p^j)$
denote the corresponding solutions to
(\ref{rid:fem2}), \eqref{Neq:sens} for $p^j$.  We define the snapshot set \blue{$Y\in\R^{n\times (m n_{max}+1})$}
\[
\blue{Y}:=\mbox{span}\left\{\left\{{\bf u}_h(p^j),({\bf u}^n_{h,\mu}(p^j))_{1\le
|\mu|=n\le n_{max}}\right\}_{1\le j\le m}\right\}
\]
including the states and the sensitivities up to order $n_{max}$
and determine a POD \blue{reduced space}
$\mathcal{V}^\ell:=\mbox{span}\{\psi_1,\ldots,\psi_\ell\}$ of rank $\ell$ by
solving the following minimization problem
\begin{eqnarray}\label{prb:pod}
\min_{\psi_1,\ldots,\psi_\ell} 
\sum_{n=0}^{n_{max}} \sum_{|\mu|=n}
\sum_{j=1}^m \beta_j \left\|{\bf u}^n_{h,\mu}(p^j)-
\sum_{i=1}^\ell \langle {\bf u}^n_{h,\mu}(p^j),\psi_i\rangle_{\bf W} \psi_i\right\|_{\bf W}^2\nonumber\\
\mbox{ s.t.} \langle\psi_j,\psi_i\rangle_{\bf W}=\delta_{ij}\quad\mbox{for }1\leq i,j\leq \ell,
\end{eqnarray}
where we set
${\bf u}^0_{h,\mu}(p^j)={\bf u}_h(p)$ with $|\mu|=0$ for brevity,
$n_{max}$ is the maximum order of considered sensitivies, $\beta_j$ are nonnegative weights, $\delta_{ij}$ denotes the Kronecker symbol, ${\bf W}$ is a  symmetric positive definite $N\times N$ matrix and $\psi_i\in\mathbb R^N$. The weighted inner product used  is defined as follows: $\blue{\langle {\bf u}, {\bf v} \rangle_{\bf W}} := {\bf u}^\top {\bf W} {\bf
v}$.

It is well-known (see \cite{GV13}) that problem (\ref{prb:pod}) admits a unique solution $\{\psi_1,\ldots,\psi_\ell\}$, where $\psi_i$ denotes the $i-$th eigenvector of the self-adjoint linear operator $\mathcal R:\mathbb R^n \rightarrow \mathbb R^n,$ i.e., $\mathcal{R}\psi_i=\lambda_i\psi_i$ with $\lambda_i\in\R$ non-negative, where $\mathcal R$ is defined as
$$
\mathcal{R}\psi=\sum_{k=1}^{n_{max}}\sum_{|\mu|=n}\sum_{j=1}^m
\beta_j \langle {\bf u}^n_{h,\mu}(p^j),\psi\rangle_{\bf W}
{\bf u}^n_{h,\mu}(p^j) \quad \mbox{for } \psi\in \mathbb R^n.
$$
Furthermore, the error in \eqref{prb:pod} can be expressed as
\begin{equation}\label{err:svd}
\sum_{n=0}^{n_{max}} \sum_{|\mu|=n}
\sum_{j=1}^m \beta_j \left\| {\bf u}^n_{h,\mu}(p^j)-
\sum_{i=1}^\ell \langle {\bf u}^n_{h,\mu}(p^j),\psi_i\rangle_{\bf W} \psi_i\right\|_{\bf W}^2
=\sum_{i=\ell+1}^d \lambda_i,
\end{equation}
where $d$ is the rank of the snapshots matrix \blue{$Y$}.

\subsection{POD approximation for state and sensitivities}

We briefly recall how to generate the reduced order modeling by means of POD. Suppose we have computed the POD basis $\{\psi_i,\ldots,\psi_\ell\}$ of rank $\ell$ according to the minimization problem \eqref{prb:pod}. For the weight matrix we choose
$$
{\bf W}: = {\bf K}(\bar p) + {\bf M}(\bar p),
$$
where $\bar p$ is a fixed reference parameter and $\bf M$ denotes the mass matrix. Then, ${\bf W}$ is the matrix associated to the discrete $H^1-$norm. We define the POD ansatz for the state as ${\bf u}_h^\ell(p):=\sum_{i=1}^\ell ({\bf \bar u}^\ell)_i\psi_i.$ This ansatz in (\ref{rid:fem2}) leads to an $\ell-$dimensional linear system for the
unknowns $\{{\bf \bar u}_i\}_{i=1}^\ell$, namely
\begin{equation}\label{disc:pod}
{\bf K^\ell}(p){\bf \bar u}(p)={\bf f^\ell}(p).
\end{equation}
Here the entries of the stiffness matrix are given by ${(\bf K^\ell})_{ij} = \psi_j^\top{\bf K}(p)\psi_i$ for $1\le i,j \le\ell$. The right hand side is given by ${(\bf f^\ell})_{i} = \psi_i^\top{\bf f}(p)$, $1\le i \le \ell$. We recall that due to the previously introduced  affine decomposition this projection has to be computed only once and the system matrix can be written as
$$
{\bf K}^\ell(p) = \sum_{q = 1}^{Q} \hat\Phi_q^a(p) {\bf K}^{q,\ell},
$$
where $({\bf K}^{q,\ell})_{ij} = \psi_i^\top {\bf K}^q\psi_j$ for $1 \le i,j \le \ell$ and $q=1,\ldots,Q$. The same structure can be used for the right hand side. Note that this is very important in order to obtain an efficient reduced order
model, since the system can be set up for different values of $p$ without the need of the original high dimensional matrices and right hand sides. The reduced order state equation reads:
\begin{equation}\label{rid:pod}
\left\{\begin{array}{l}
\hbox{For }p\in\mathcal{D}\hbox{ find } {\bf u}_h^\ell\in \blue{\R^\ell}\hbox{ s.t.}\\
{\bf K}^\ell(p){\bf u}_h^\ell(p)={\bf f}^\ell(p).
\end{array}\right.
\end{equation}
In an analogous way we obtain the general reduced sensitivity equation from
(\ref{Neq:sens}).  We need to make an ansatz for the sensitivities ${\bf
u}^{n,\ell}_{h,\mu}$ and project the system onto the subspace spanned by the
POD basis.  In the present work we use the same basis
functions for the state and the sensitivity variables and to achieve a
better approximation property of the reduced order sensitivities
we add the solution of the sensitivity
equation to the snapshots set. Note that then the stiffness matrix
in the reduced sensitivity
equation is the same as in the reduced state equation,
so that only the right-hand side needs to be projected.

%Finally we note that the reduced sensitivity equation; due to the structure of the problem the stiffness matrix is the same of the reduced state equation \eqref{disc:pod}, and only the right hand-side needs to be projected.

\subsection{A-posteriori error estimations}

A-posteriori error estimators have a crucial role in model order reduction. They provide a certification of the surrogate model without the need of the computation of the truth solution. In the present work, we consider as truth solution the finite element approximation. 
If the transformation on a subdomain is given by \eqref{exp_dec}, the bilinear form reads
$$a(u,v;p)=\sum_{q=1}^Q\sum_{i,j=1}^2 [\left(C_q(p)\right)^{-1}\nu_q \left(C_q(p)\right)^{-T}]_{ij} |\blue{\det C_q(p)}| \int_{\Omega_q(\bar p)} \dfrac{\partial u}{\partial x_i}\dfrac{\partial v}{\partial x_j} dx,$$
and a lower bound for the coercivity constant $\alpha(p)$ in \eqref{coer} can be computed from the following estimate
\begin{equation}\label{est_coer}
a(v,v;p) \ge \min_q\left\{|\det(C_q(p))| \blue{\lambda_{min}}(C_q(p)^{-1} \nu_q C_q(p)^{-T})\right\} a(v,v;\bar p), 
\end{equation}
where $\blue{\lambda_{min}}$ denotes the smallest eigenvalue of the matrix $\left(C_q(p)^{-1}\nu_q C_q(p)^{-T}\right)$ 
%since
%$$\nabla \hat v^T (C_q(p)^{-1} C_q(p)^{-T}) \nabla \hat v  \geq a_{min}(C_q(p)^{-1} C_q(p)^{-T}) \nabla \hat v^T \nabla \hat v^T$$
%
%
%following the
%min-$\hat\Phi$ theorem introduced in \cite{PR07}. It turns out that 
%%
%$$
%\blue{
%\alpha_{LB}(p):=\min_{\substack{q=1,\ldots,Q}} \dfrac{\hat\Phi^{q}_a(p)}{\hat\Phi^{q}_a(\bar{p})}\le\alpha(p),
%}
%$$
%
and $\bar{p}$ is the fixed reference parameter. %From the estimate \eqref{est_coer} it is possible to compute a lower bound for the coercivity constant. 
An upper bound for the continuity constant $\gamma(p)$ in \eqref{cont:cost} can be computed with the same approach.

%is given by 
%%
%$$
%\blue{
%\gamma(p)\le\gamma_{UB}(p):=\max_{\substack{q=1,\ldots,Q}} \dfrac{\hat\Phi^{q}_a(p)}{\hat\Phi^{q}_a(\bar{p})}.
%}
%$$
%
% \blue{We note that $\alpha_{LB}(p)$ and $\gamma_{UB}(p)$ are computed by taking the minimum, respectively maximum, over the two components generated by the minimum/maximum over $q$. This is a consequence of our assumption on isotropic materials.}
%For more details we refer the reader to \cite{PR07}. \red{Further, we note that under our assumptions the $\min-\hat\Phi$ theorem hold true as explained in \cite{H17}.} %We refer the reader to \cite{PR07} for more details on the approximation of the coercivity and continuity constants.
Then, we can derive an error bound for the reduced state and sensitivity equations (see \cite{RHP08} for more details) in terms of the reduced residual of the aforementioned equations. 

For this purpose we define the residual for equation \eqref{Neq:sens} as follows:

\begin{definition}
Let $u^{n,\ell}_{h,\mu}\in \blue{X_h(\bar\Omega)}$, $n = 0,\ldots,n_{max}$,
be the reduced order solution of \eqref{Neq:sens}. We define the residual
\begin{equation}\label{nsensres}
r_{u^n_\mu}(v;p):=
 \frac{\partial^\mu f}{\partial p^\mu}(v;p) 
- \sum_{1\le |\kappa|\le n} \binom{\mu}{\kappa}
\frac{\partial^\kappa a}{\partial p^\kappa}(u^{n-|\kappa|,\ell}_{h,\mu-\kappa},v;p)
-a( u^{n, \ell}_{h,\mu},v;p).
\end{equation}
\end{definition}

%We note that for $n=0$ and $n=1$ we find the state and sensitivity residual \eqref{stateres} and \eqref{sensres}. 

Then, we have
\begin{theorem}
Let $u_{h,\mu}^n\in \blue{X_h((\bar\Omega)}$ be the solution to \eqref{Neq:sens} for
$|\mu|=n$ and
$u^{n, \ell}_{h,\mu}\in \blue{X_h(\bar\Omega)}$ be the corresponding reduced solution of
\eqref{Neq:sens}. Then, the following inequality holds
\begin{equation}\label{Napost_sens}
\| u_{h,\mu}^n- u_{h,\mu}^{n, \ell}\|_{\blue{X_h(\bar\Omega)}}\leq
\Delta_{u_\mu^n}^{\ell}(p):=\dfrac{1}{\alpha(p)}
\left(\|r_{u^n_\mu}(p)\|_{\blue{(X_h(\bar\Omega)})'}
+\sum_{1\le |\kappa|\le n} \binom{\mu}{\kappa}
\gamma_\kappa\Delta^\ell_{u^{n-|k|}_{\mu-\kappa}}(p)\right),
\end{equation}
%
%where $\Delta u^\ell_n:=\dfrac{\|r_{u^n}(\cdot\,;p)\|_{(X_h)'}}{\alpha(p)}$ and 
where $\gamma_\kappa$ is the continuity constant of the
$\kappa$-th derivative of the  coercive bilinear form.
\end{theorem}

{\em Proof.} We denote the error by $e_{u^n_\mu}:=  u^n_{h,\mu} -
u^{n,\ell}_{h,\mu}$. With \eqref{eq:sens_general} we find that
\begin{eqnarray*}
a(e_{u^n_{h,\mu}},v;p)&=&a( u_{h,\mu}^n,v;p)-a( u^{n, \ell}_{h,\mu},v;p)\\
&=&\frac{\partial^\mu f}{\partial p^\mu}(v;p)
- \sum_{1\le |\kappa|\le n} \binom{\mu}{\kappa}
\frac{\partial^\kappa a}{\partial p^\kappa}(u_{h,\mu-\kappa}^{n-|\kappa|},v;p)
-  a( u^{n, \ell}_{h,\mu},v;p)\\
&=& r_{u_\mu^n}(v;p) + \sum_{1\le |\kappa|\le n} \binom{\mu}{\kappa}
\frac{\partial^\kappa a}{\partial
p^\kappa}(u_{h,\mu-\kappa}^{n-|\kappa|,\ell}-u_{h,\mu-\kappa}^{n-|\kappa|},v;p),
\end{eqnarray*}
Now we set $v = e_{u^n}$ and obtain
$$
a(e_{u^n_{h,\mu}},e_{u^n_{h,\mu}};p) = 
r_{u^n_\mu}(e_{u^n_{h,\mu}};p) + \sum_{1\le |\kappa|\le n}
\binom{\mu}{\kappa}
\frac{\partial^\kappa a}{\partial p^\kappa}(u_{h,\mu-\kappa}^{n-|\kappa|, \ell} -
u_{h,\mu-\kappa}^{n-|\kappa|},e_{u^n_{h,\mu}};p) 
$$
By applying Cauchy-Schwarz and using the coercivity of $a$ as well as the
continuity of $\frac{\partial^\kappa a}{\partial p^\kappa}$ we find
\begin{align*}
\alpha(p)\|e_{u^n_{h,\mu}}\|^2_{\blue{X_h(\bar\Omega)}} &\le
\|r_{u^n_\mu}(p)\|_{\blue{(X_h(\bar\Omega))}'}\|e_{u^n_{h,\mu}}\|_{\blue{X_h(\bar\Omega)}}\\
&+ \sum_{1\le|\kappa|\le n} \binom{\mu}{\kappa}
\gamma_\kappa\|u_{h,\mu-\kappa}^{n-|\kappa|, \ell} -
u_{h,\mu-\kappa}^{n-|\kappa|}\|_{\blue{X_h(\bar\Omega)}}\|e_{u^n_{h,\mu}}\|_{\blue{X_h(\bar\Omega)}}.
\end{align*}
Dividing by $\alpha(p)\|e_{u^n_{h,\mu}}\|_{\blue{X_h(\bar\Omega)}}$ leads to
\eqref{Napost_sens} first for $n=0$ and by inductively using the bound
$\|u_{h,\mu-\kappa}^{n-|\kappa|, \ell} - u_{h,\mu-\kappa}^{n-|\kappa|}\|_{\blue{X_h(\bar\Omega)}}\le
\Delta^\ell_{u_{\mu-\kappa}^{n-|\kappa|}}$
for $n\le n_{max}$. $\Box$

\begin{rmk} 
We note that \eqref{Napost_sens} is a generalization of well-known error bounds for state and first order sensitivity equations, compare e.g. \cite{OP07}, and \cite{DH15} for time-dependent problems.
\end{rmk}
\begin{rmk}
%The presented generalized error estimator \eqref{Napost_sens} is not limited to the introduced setting arising from a parametrized shape optimization. In fact, 
If only the linear form $f$ depends on the parameter $p$, i.e., $a(u,v) = f(v;p)$, the error estimator \eqref{Napost_sens} reads
$$
\| u_{h,\mu}^n- u_{h,\mu}^{n, \ell}\|_{\blue{X_h(\bar\Omega)}}
\leq\dfrac{\|r_{u_\mu^n}(p)\|_{(\blue{X_h(\bar\Omega)})'}}{\alpha}.
$$
\end{rmk}

\blue{We note that the solution $u(\mu)$ in this case lives in a $Q_f$-
dimensional linear subspace, where $Q_f$ is defined in \eqref{source_aff}, see e.g. \cite{H17}. Then, for $\ell\ge Q_f$ the error with the reduced model is zero and no error bound is needed.} We note that sensitivity equations and error estimators for the uncertainties $\phi$ are analogous to \eqref{nsensres}-\eqref{Napost_sens} and are left to the reader.

\subsection{The POD method for optimization problem}
In this section we explain how to solve the parametrized optimization problem.
In our construction the POD spaces depend on the parameter $p$. Since the solution to the parameter optimization problem is not known in advance the POD space has to be adapted/enriched during the parameter optimization procedure. To achieve this goal we here propose a certified extension of the approach suggested in \cite{AH01} to our robust setting, where the error bound, $\Delta_{u_\mu^n}^{\ell}(p)$ introduced in Section 4.2, helps in the selection of the snapshot sets.
%To build the reduced order model we have to select some relevant parameters related to \red{optimization problem.} The snapshot set is important for POD model reduction, since the basis functions are built upon this set. Therefore, the choice of the parameters should be done according to our \red{optimization problem.} For this reason we propose a goal-oriented algorithm for the computation of the POD basis functions \new{where the goal is the solution of the optimization problem. A similar approach for time-dependent PDE constrained optimization problem can be found in \cite{AH01} where it was proposed an adaptive method for POD but they did not provide any certification for the problem considered.}
%\new{The error bound $\Delta^\ell_{u^n}$ for state and sensitivity introduced in Section 4.2 helps also in the selection of the snapshots set.}
The algorithm works as follows: we start with a very coarse parameter sample choosing only one parameter $p^0$ and solve the full problem together with the sensitivity equations associated to this parameter. Then, we compute the POD basis functions and perform the reduced optimization procedure. At the end of the process we find a new parameter $p^1$ which is an approximation of the optimal desired design, we update the parameter set $\mathcal{D}=\{p^0,p^1\}$, solve the full problem and the sensitivity equations related to the new parameter $p^1$. Then, we enlarge the snapshots set and compute new POD basis functions. We iterate this process until the stopping \blue{criterion} is reached. The procedure is summarized in Algorithm \ref{Alg:POD}.

\begin{algorithm}[H]
\caption{(Adaptive POD optimization)}
\label{Alg:POD}
\begin{algorithmic}[1]
\REQUIRE $p^0, tol>0,$
\STATE Set $k =0, \mathcal{V} = [\,],$
\STATE Set Snapshot set 
$$\mathcal{V}=\left[\mathcal{V}, \left\{{\bf u}_h(p^k),({\bf u}^n_{h,\mu}(p^k))_{1\le
|\mu|=n\le n_{max}}\right\}\right],$$
\STATE Compute POD basis functions $\{\psi_i\}_{i=1}^\ell$ with $\ell=\mbox{rank}(\mathcal{V})$
\STATE Find $p^{k+1}$ solving the OCP with the reduced order model \eqref{rid:pod}
\IF{$\max\limits_{|\mu|\le n_{max}} \left(\Delta_{u_\mu^n}^{\ell}(p^k)\right) > tol$}
\STATE Set k=k+1
\STATE GOTO 2
\ENDIF
%\STATE Compute ${\bf u}^{i,\ell}(p^{k+1})$ and ${\bf u}^i(p^{k+1})$, 
%\ENDWHILE
\end{algorithmic}
\end{algorithm}

%
%\begin{algorithm}[H]
%\caption{(Goal-Oriented POD optimization)}
%\label{Alg:POD}
%\begin{algorithmic}[1]
%\REQUIRE $p^0, {\bf u}^0(p^0),\mathcal{V}=[], k=0, tol>0$
%\WHILE{\red{$\Delta^\ell_{u^n}(p^k)\geq tol$ for $n=0,1,\ldots, n_{max}$}}
%\STATE Compute sensitivities ${\bf u}^i(p^k)$, for $i=0,1,2,3,\ldots \blue{n_{max}}$
%\STATE Set Snapshot set 
%$$\mathcal{V}=[\mathcal{V}, {\bf u}^0(p^k), {\bf u}^1(p^k), {\bf u}^2(p^k), {\bf u}^3(p^k),\ldots {\bf u}^{N_p}(p^k)]$$
%\STATE Compute POD basis functions $\{\psi_i\}_{i=1}^\ell$ with $\ell=\mbox{rank}(\mathcal{V})$
%\STATE Find $p^{k+1}$ solving the OCP with the reduced order model \eqref{rid:pod}
%\STATE Compute ${\bf u}^{i,\ell}(p^{k+1})$ and ${\bf u}^i(p^{k+1})$, 
%\STATE Set k=k+1
%\ENDWHILE
%\end{algorithmic}
%\end{algorithm}
%
In our simulations this approach turned out to be very efficient since it avoids long pre-computations. In this way we are able to update the snapshot set and the POD basis functions. Our update contains information on the \red{optimization problem} and it improves the quality of our surrogate model.
Note that every reduced optimization problem contains the a-posteriori error for the state and sensitivity equations introduced in Section 4.2. %In this way we update the snapshot set, and therefore the POD basis, not only when we are close to the minimum of the problem, but even if the reduced model is not accurate enough. %However, our update contains information on the optimal control problem and it improves the quality of our surrogate model. 
%\myc{MH5: Why is it called 'goal oriented'? Which is the goal? Imo this is APOD as in \cite{AH01}, We then clearly should state that the goal is the solution of the optimization problem. In the original work apod adds snapshots to the existing snapshot sets, see 135 in my publication list. Then the presented algorithm is apod for the robust optimization problem. }

% -------------------------------------------------
% -------------------------------------------------
% --------------- SECTION 5 -----------------
% -------------------------------------------------
% -------------------------------------------------
\newpage
\section{Numerical tests}\label{sec:test}

\begin{figure}[htbp]
\centering
\includegraphics[width=5.8cm]{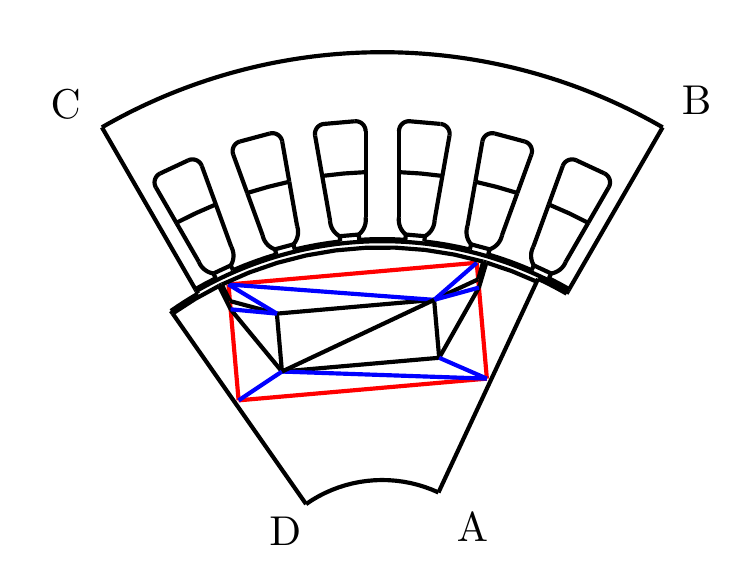}
\includegraphics[width=6cm]{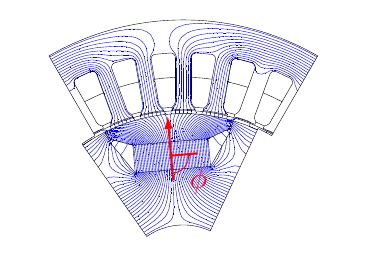}
\caption{\label{fig:geometry}Geometry configuration with region for the affine decomposition (red box) and triangulation for the decomposition (blue lines) (left plot). Magnetic vector potential for the geometry configuration and magnetic field angle $\phi$ (right plot).}
\end{figure}

In our numerical example we consider an optimal design problem for a permanent magnet synchronous machine. We start by introducing the model and geometry under consideration. We consider a three-phase six-pole permanent magnet  synchronous machine (PMSM) with one buried permanent magnet per pole. The geometry is shown in Figure~\ref{fig:geometry}. The goal of the design optimization is to change the size and location of the permanent magnet such that the material of the magnet is minimized while maintaining \blue{a lower bound on} the electromotive force. We consider a description using three parameters: $p_1$ the width, $p_2$ the height and $p_3$ the central perpendicular distance between the rotor and the surface of the magnet. The region around the permanent magnet (Figure~\ref{fig:geometry} red box) is decomposed into twelve triangles (Figure~\ref{fig:geometry} blue lines). The introduced triangulation of the parametrized domain allows to perform the affine linear decomposition as introduced in \eqref{eq:affine} with $Q_a=12=Q_f$.

PMSMs can be described sufficiently accurate by the magneto-static approximation of Maxwell's equation. \blue{The governing} parametrized equation is given by 
\begin{equation}
\label{eq:model_problem}
\nabla\times\left(\nu \nabla\times u(p)\right) = J_{src}(p)  - \nabla\times H_{pm}(p), 
\end{equation}
with boundary conditions
$$
u|_{BC} = u|_{DA} = 0\quad\mbox{and}\quad u|_{AB} = -u|_{CD}
$$
where $\nu$ is the reluctivity, $J_{src}$ is the source current density and $H_{pm}$ the field of the permanent magnets (PM). We note that equation \eqref{eq:model_problem} fits into the abstract formulation \eqref{pde_par}. In the 2D planar setting together with a finite element method for the magnetic vector potential, we lead to the discrete form given by the linear model presented in \eqref{rid:fem2}. %where the bilinear form $a$ and the right hand side $f$ are given as \myc{'fill in the forms here'}. 
To extract performance values the loading method is used to exploit the frequency domain \cite{RZ91}. To obtain quantities like the electromotive force (EMF) a Fourier analysis of the magnetic vector potential around the inner surface of the stator is carried out. This can be written as a linear function $E_0  =  \mathbb E^\top {\bf u}_h$. More details on the configuration we adopt in the present work can be found in \cite{AHLU15,HPHB97,P98}. 

Let us next formulate the optimization problem. We start by introducing the nominal optimization problem. The goal of the optimization is to minimize the required material for the permanent magnet while maintaining \blue{a lower bound on} the EMF $E_0$. In the mathematical model this leads to a cost function of the form
$$
\min_{p\in\mathbb R^3} g_0(p) := p_1 p_2 + \rho \max(0,E_0^d - E_0(p,{\bf u}(p))),
$$
where $E_0^d$ is the desired \blue{lower bound on the} EMF and $\rho\in\R^+$ a weight parameter. Additionally, we have the constraints
$$
(1,1,5) \le (p_1,p_2,p_3) \le (\infty,\infty,14),\quad p_2 + p_3 \le 15 \quad\mbox{and}\quad 3p_1 - 2p_3 \le 50.
$$
The upper and lower bounds for the parameters and the first inequality are due to the parametrization of the geometry and the restriction that the permanent magnet has to stay within the red box in Figure~\ref{fig:geometry} (left). The last inequality is a design restriction that avoids that the corner of the permanent magnet comes too close to the rotor surface. Note that we right away use the reduced formulation as introduced in \eqref{eq:opt_red}. To obtain a smooth formulation we introduce a slack variable. We reformulate the nominal optimization problem by using the variable $x=(p,\xi)$ as
\begin{equation}
\label{eq:opt_nom_ex}
\min_{x\in\mathbb R^4} g_0(x) := p_1 p_2 + \rho \xi
\quad\mbox{s.t.}\quad 
g_{1,\ldots,8}(x)=\left(\begin{array}{c}
       p_2 + p_3 - 15\\
       3p_1-2p_3-50\\
       E_0^d - E_0({\mathbf u_h},p) - \xi\\
       1 - p_1\\
       1 - p_2\\
       5 - p_3\\
       -\xi\\
       p_3 - 14
      \end{array}\right) \le 0.
\end{equation}
In this form the optimization problem fits exactly into the framework of \eqref{eq:opt_red}.

Next let us introduce the uncertainty. For our numerical example we assume that the magnetic field  angle in the permanent magnet is uncertain \cite{OH13}. This can be due to manufacturing imprecision. In the nominal optimization the magnetic field is aligned perfectly, i.e. the field angle is $90$, see Figure~\ref{fig:geometry} (right plot). In practice this can not be met and a deviation is to be expected. The field angle enters the model problem \eqref{eq:model_problem} nonlinearly through the right hand side, in particular in the term $H_{pm}$. In our numerical example we allow a field angle $\phi$ in the range \new{$[85, 95]$}. Following the definition of the uncertainty set in Section~\ref{sec:opt_robust} we get
$$
\mathcal U_k = \left\{\phi\in \mathbb R\,\big|\,\phi = \new{90} + \delta,\,\|0.2\delta\|_k \le 1  \right\}
$$
with $k\in\{2,\infty\}$, i.e. $\hat\phi = \new{90}$ and $D = 5$ in \eqref{un_set}. Using \blue{these} settings we can now solve the linear and quadratic approximation of the robust counterpart \eqref{eq:robust_approx_linear} and \eqref{eq:robust_approx_quadratic}. In the case of the linear approximation we choose $k = \infty$ and for the quadratic approximation we set $k = 2$. Note that in the presented setting these two uncertainties sets are identical. The different choices are for technical reason as outlined in the derivation.

Before presenting the numerical results let us give a short overview of the numerical strategy utilized to solve the optimization problems. The computations are carried out in MATLAB. To solve the nominal and robust optimization problems an SQP method with Armijo-backtracking strategy using a $\ell_1$-penalty function is used \cite{NW06}. The Hessian is computed via BFGS updates. Alternatively also routines like {\tt fmincon} in MATLAB can be used obtaining similar results. The derivative of $g_i,\,i=0,\ldots,8$ are computed using the sensitivity approach \cite{HPUU09}. Also the derivatives with respect to the uncertain parameter $\phi$ are computed using this approach. The structure of the sensitivity equations are as outlined in \eqref{eq:sens}. 

\subsection{Results obtained by the finite element approximation}

We start by presenting the numerical results utilizing the finite element approximation. Piecewise linear and continuous finite elements are used to discretize equation \eqref{eq:model_problem} leading to a system with $61013$ degrees of freedom.

The initial geometry configuration corresponding to $p = (19,7,7)$ is shown in Figure~\ref{fig:geometry} (left) together with the corresponding magnetic vector potential (right). From the magnetic vector potential we extract the EMF which we will use as the desired value $E_0^d = 30.34$ in our optimization problem. 

In Table~\ref{tab:results1} we show the results obtained in the optimization. In column $V_{pm}$ the volume of the permanent magnets is given. Note that through the optimization a significant reduction in size is achieved. The ratio is given in percent in the second column. For the nominal optimization a reduction by $53\%$ is obtained and in the robust case reductions of $50\%$ (linear case) and $49\%$ (quadratic case). The respective parameters are given in the third column. Lastly, in column four and five the EMF $E_0$ for the field angle $\phi = 90$ and the worst-case EMF $E_0^{worst}$ are given. It can be seen that the uncertainty in the magnetic field of the permanent magnet has an impact on the performance. In the robust optimization this influence is incorporated. Hence it can be seen that in the case of the quadratic approximation \new{(Rob.\ Quad and Rob.\ Shift)} very good results can be achieved: the worst-case EMF stays above the target value of $30.34$.
\new{Shifting the expansion point (Rob.\ Shift) helps to avoid being overly conservative, so that the EMF target is precisely met, which allows to reduce the volume slightly.} In the case of the nominal optimization a significant decay can be observed. In the case of the linear approximation the target can not be met since the approximation is not accurate enough. It can also be observed that by performing only a nominal optimization the worst-case can decrease compared to the initial configuration. In the presented case the difference is small but can become more significant in different settings.

\begin{table}
\centering
\begin{tabular}{lccccc}
 \toprule
            & $V_{\mathrm{pm}}$ & \%  & $p$ & $E_0$ ($90$) & $E_0^{\mbox{worst}}$\\
 \midrule
 Init       & $133.00$           & $100$           & $(19.00,7.00,7.00)$    & $30.3483$ & $29.8873$\\
 Nom. Opt   & $\phantom{1}62.36$ & \phantom{1}$47$ & $(21.08,2.96,6.62)$    & $30.3483$ & $29.8873$\\
 Rob. Lin   & $\phantom{1}66.46$ & \phantom{1}$50$ & $(21.10,3.15,6.64)$    & $30.5817$ & $30.2322$\\
 Rob. Quad  & $\phantom{1}67.95$ & \phantom{1}$51$ & $(21.10,3.22,6.65)$    & $30.6995$ & $30.3487$\\
 \new{Rob. Shift} & \new{$\phantom{1}67.94$} & \new{\phantom{1}$51$} & \new{$(21.10,3.22,6.65)$}    & \new{$30.6992$} & \new{$30.3483$}\\
 \bottomrule
\end{tabular}
\caption{Comparison of the results obtained by the optimization and robust optimization using different approximation orders.}
\label{tab:results1}
\end{table}

\begin{figure}[htbp]
\centering
\includegraphics[width=3.9cm]{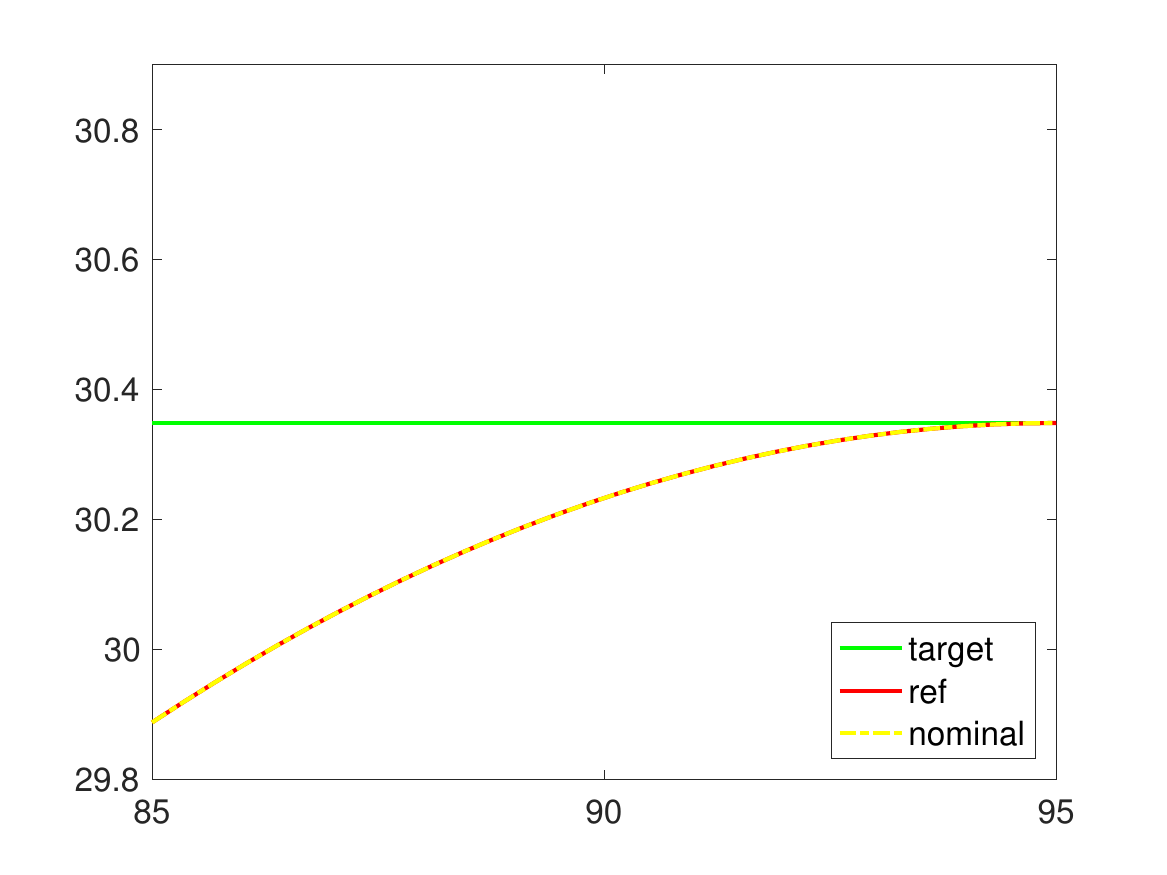}
\includegraphics[width=3.9cm]{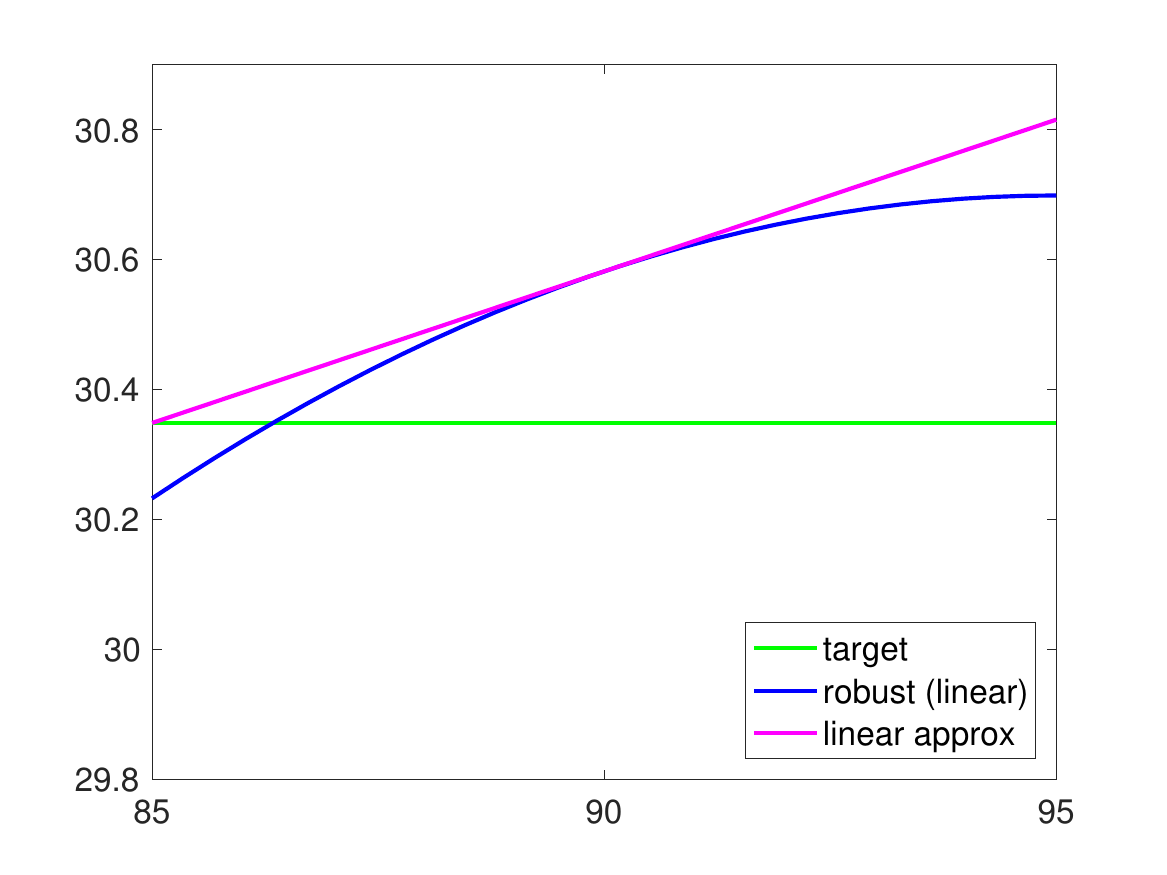}
\includegraphics[width=3.9cm]{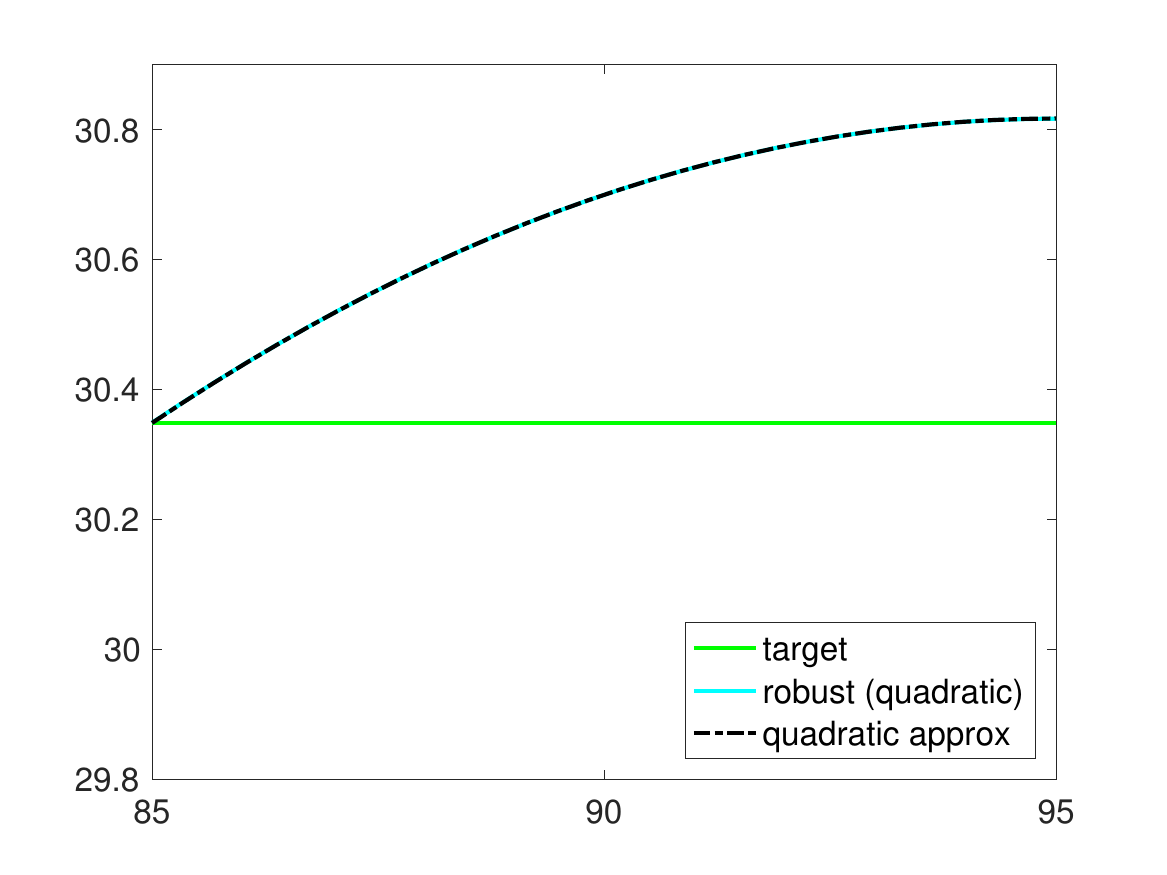}
\caption{\label{fig:emf}Influence of the permanent magnet field angle on the electromotive force for different geometry configuration. }
\end{figure}

The approximation quality and the behavior of the EMF is shown in Figure~\ref{fig:emf} for the different optimal designs. In the left plot the behavior of the EMF for different permanent magnet field angles $\phi$ is shown for the initial and the nominal optimal configuration. It can be seen that the target $E_0^d$ is only reached in $\phi = 90$. In the middle plot the linear approximation and the actual EMF are compared. It can be seen that the approximation is not accurate enough to determine the worst-case. The right plot corresponds to the quadratic approximation. It can be seen that the EMF is approximated very well by the quadratic model. Hence also good results in the robust optimization using this approximation can be expected.
\new{We note that we do not include a plot for the quadratic approximation with a shifted expansion point because the difference is not visible from the standard quadratic approximation.}

%\begin{rmk}
%Note that our choice of the field angle $\phi$ in the range $[80^\circ, 90^\circ]$ has technical reasons. The EMF has an almost symmetric behavior around the field angle of $90^\circ$. Hence the linear approximation of the robust counterpart at around $\phi = 90^\circ$ will fail since the \blue{derivatives} of $g_i,\,i=0,\ldots,8$, with respect to $\phi$ are almost zero. The quadratic approximation on the other hand has no problem in achieving a good approximation.
%\end{rmk}

\begin{figure}[htbp]
\centering
\begin{tabular}{ccc}
    \begin{tikzpicture}
        \node[anchor=south west,inner sep=0] (image) at (0,0) {\includegraphics[width=3.9cm]{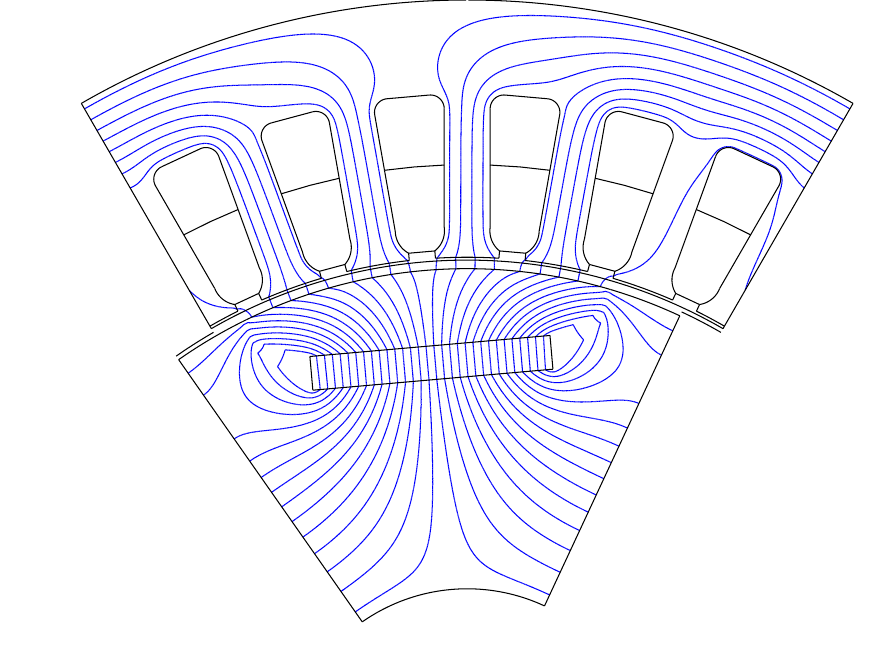}};
        %\draw[red,ultra thick,rounded corners] (80mm,35mm) rectangle (10mm,10mm);
        \begin{scope}[x={(image.south east)},y={(image.north west)}]
            \draw[red,thick] (0.5,0.3) rectangle (0.75,0.55);
            %\draw[help lines,xstep=.1,ystep=.1] (0,0) grid (1,1);
            %\foreach \x in {0,1,...,9} { \node [anchor=north] at (\x/10,0) {0.\x}; }
            %\foreach \y in {0,1,...,9} { \node [anchor=east] at (0,\y/10) {0.\y}; }
        \end{scope}
    \end{tikzpicture}
    &
    \includegraphics[width=3.9cm]{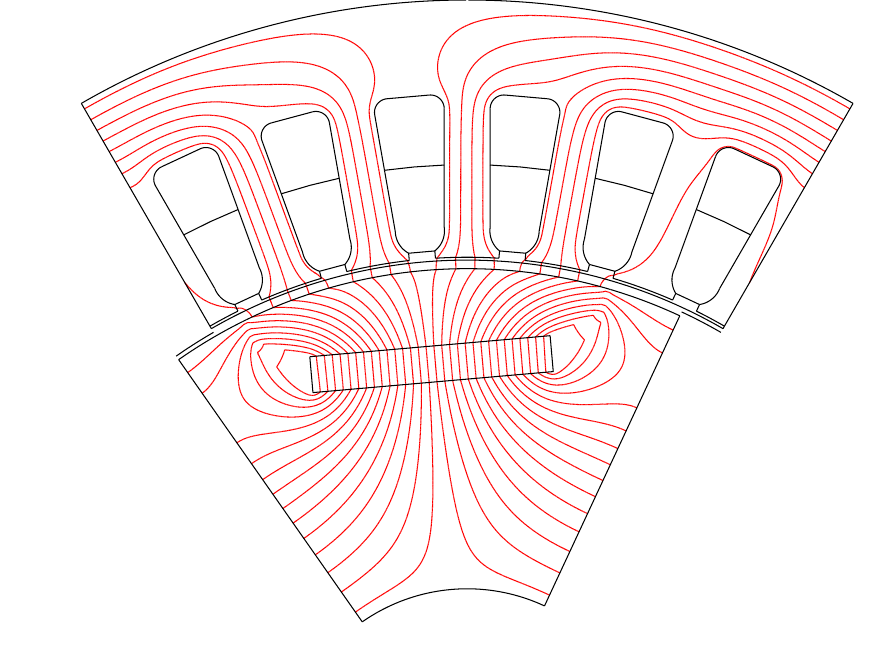}
    &
    \includegraphics[width=3.9cm]{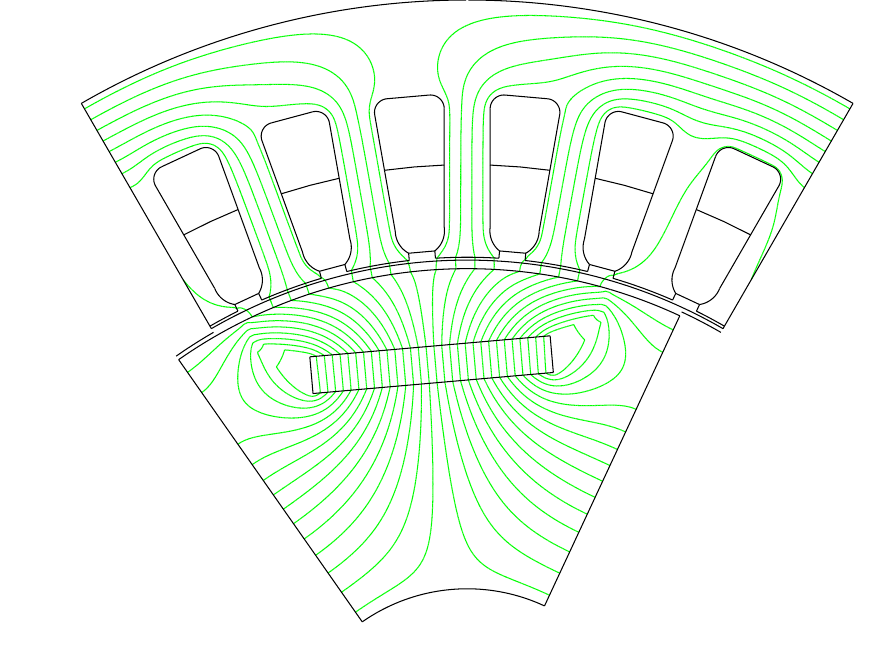}
    \\
    \scalebox{3.5}{\begin{tikzpicture}
        \begin{scope}[x={(image.south east)},y={(image.north west)}]
            \clip (0.5,0.3) rectangle (0.75,0.55);
        \node[anchor=south west,inner sep=0] (image2) at (0,0) {\includegraphics[width=3.9cm]{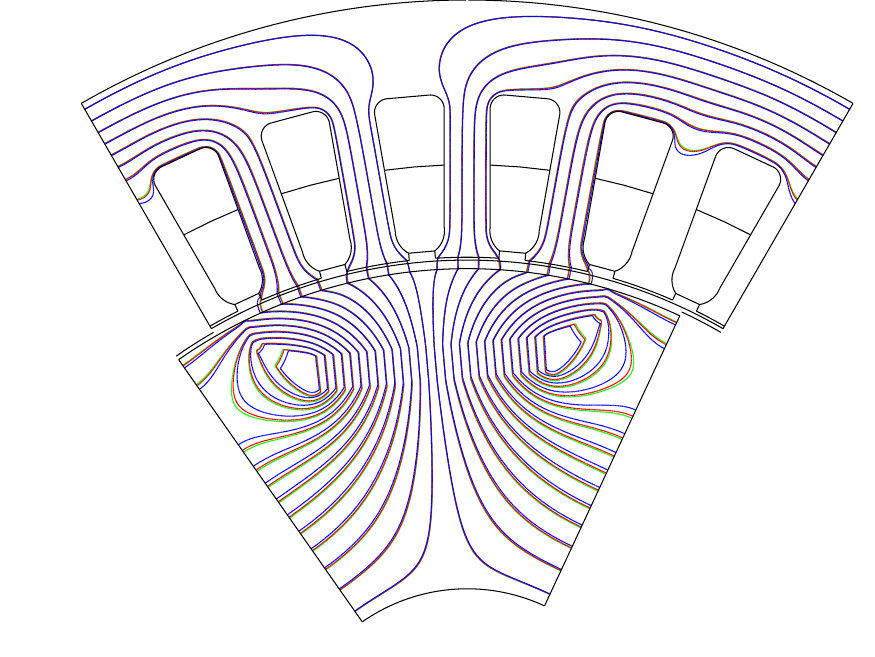}};
            %\draw[help lines,xstep=.1,ystep=.1] (0,0) grid (1,1);
            %\foreach \x in {0,1,...,9} { \node [anchor=north] at (\x/10,0) {0.\x}; }
            %\foreach \y in {0,1,...,9} { \node [anchor=east] at (0,\y/10) {0.\y}; }
            \draw [thick,red] (current bounding box.north east) rectangle (current bounding box.south west);
        \end{scope}
    \end{tikzpicture}}
    &
    \includegraphics[width=3.9cm]{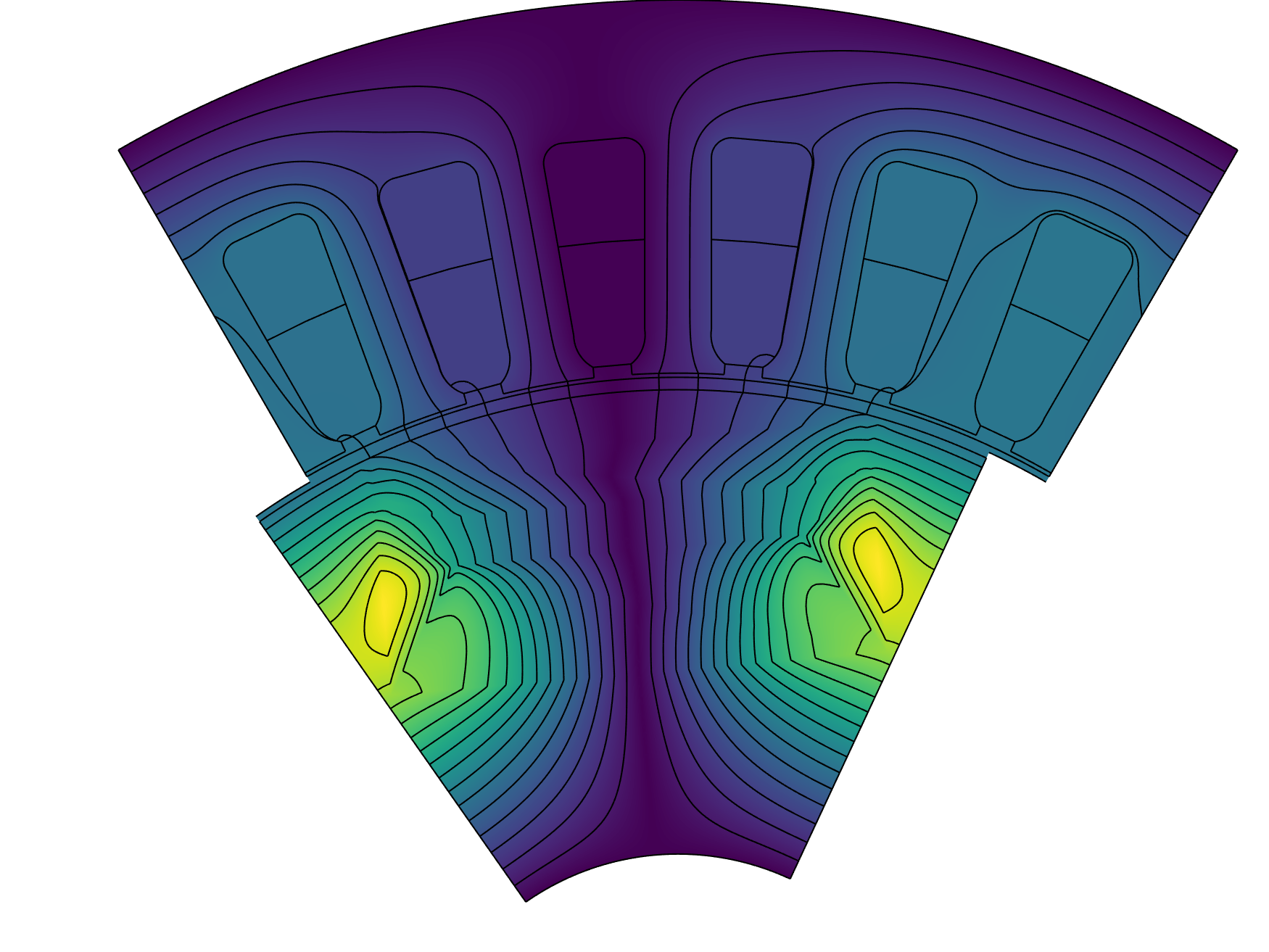}
    &
    \includegraphics[width=3.9cm]{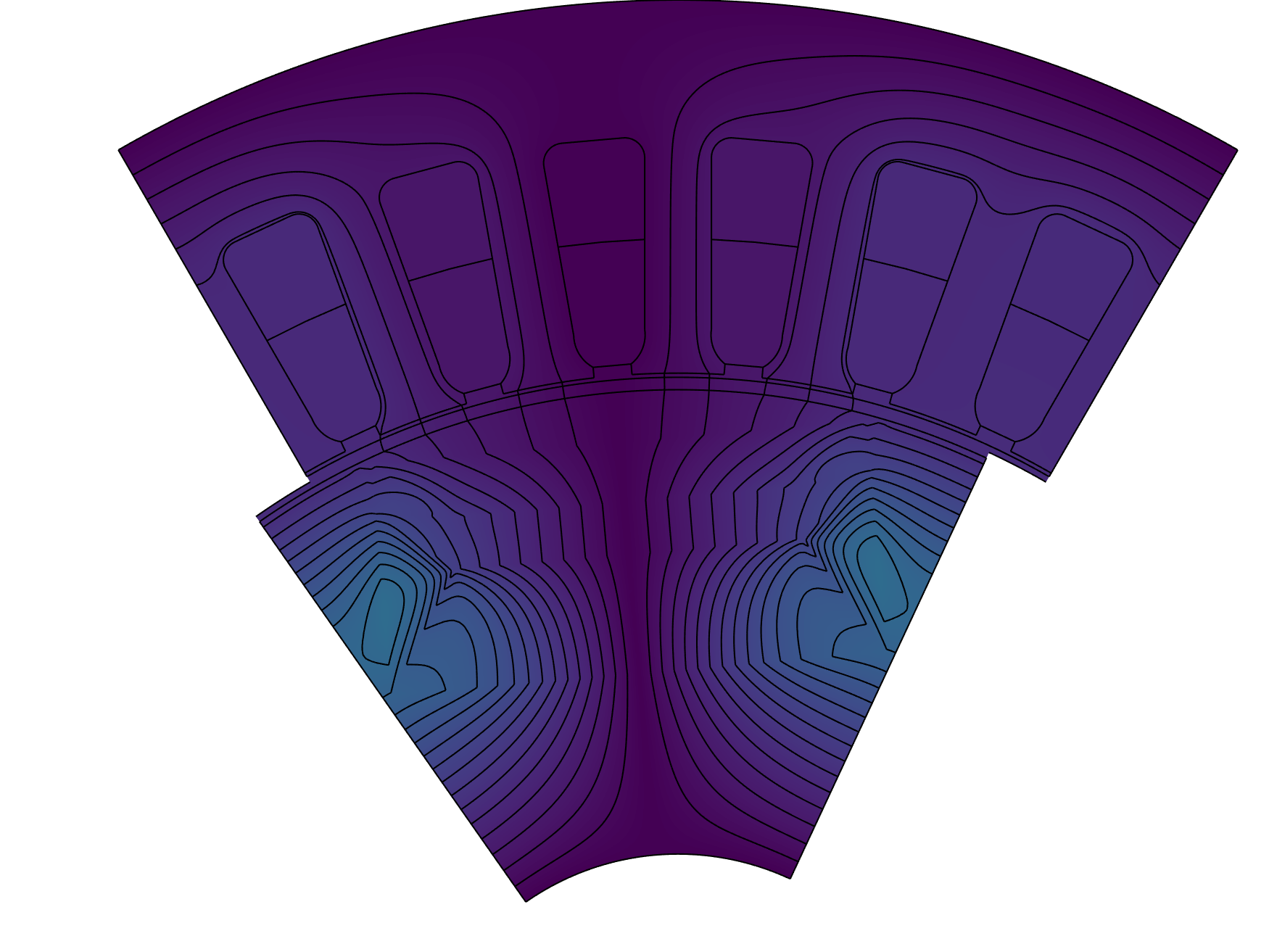}
\end{tabular}
\caption{\label{fig:mvp}Magnetic vector potential for the geometry configuration obtained by the optimization.
\blue{First row, left to right: nominal optimal design, robust optimal design with linear approximation, and with quadratic approximation.}
\blue{Second row, left to right: magnified view of the combined potentials, difference between linear and nominal, difference between quadratic and linear (brighter means larger difference).}
}
\end{figure}

In Figure~\ref{fig:mvp} the magnetic vector potentials are shown for the three optimal designs. To conclude we have a look at the computational expenses for the optimization. These are summarized in Table~\ref{tab:results2}, where the computational time in seconds, the number of iterations and the number of required PDE solves are compared. While the computational time in this example is still very low it can be seen that the robust optimization is more expensive. Especially the number of PDE solves increases significantly compared to the nominal optimization. Hence we will use model order reduction to reduce the computational costs. This will be outlined in the next section.

\begin{table}
\centering
\begin{tabular}{lccc}
 \toprule
            & CPU (s) &  Iter & PDE solves\\
 \midrule
 Nom. Opt   & $ 38.76 $  & $18$  & $72$\\
 Rob. Lin   & $42.71$    &  $8$  & $64$\\
 Rob. Quad  & $46.47$    & $6$  & $72$\\
 \new{Rob. Shift} & \new{$88.51$} &  \new{$7$} & \new{$84$}\\
 \bottomrule
\end{tabular}
\caption{Performance of the SQP method and computational cost.}
\label{tab:results2}
\end{table}
  
\subsection{Results obtained by the reduced order model}

Let us start our analysis of the reduced order model by dealing with the approximation of the state and the sensitivity variables. To this aim we select $125$ parameters chosen as follows:
$$\mathcal{D}_{train}=  \{1, 5.75,  10.5,  15.25,  20\}\times \{1, 2, 3, 4, 5\} \times \{4, 5.5, 7, 8.5, 10 \} $$
and we compute an enlarged snapshot set with state and sensitivities for each parameter in $\mathcal{D}_{train}$. Altogether, we collect $500$ snapshots. We note that we do not need to compute a reduced order model with respect to the slack variable $\xi$ since it is not a variable that is involved in the PDE.

\begin{figure}[htbp]
\centering
\includegraphics[width=6cm]{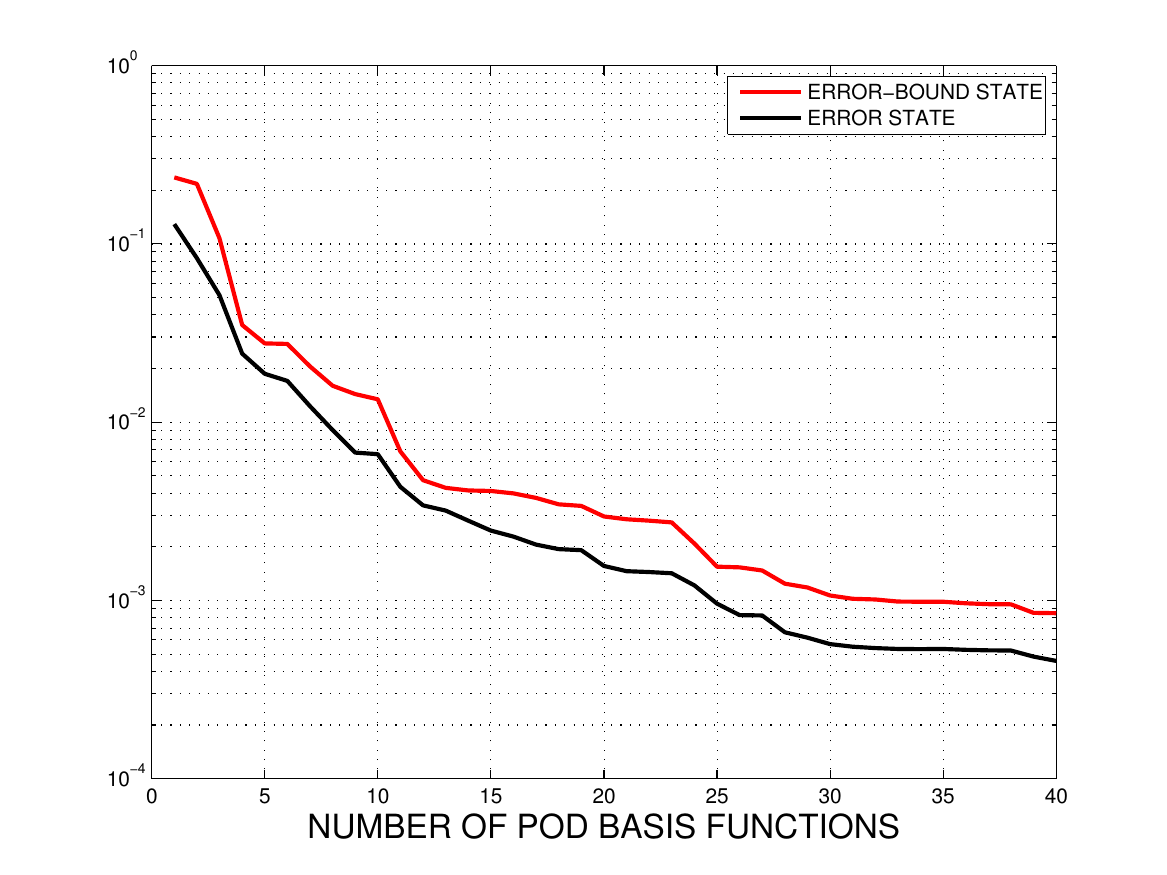}
\includegraphics[width=6cm]{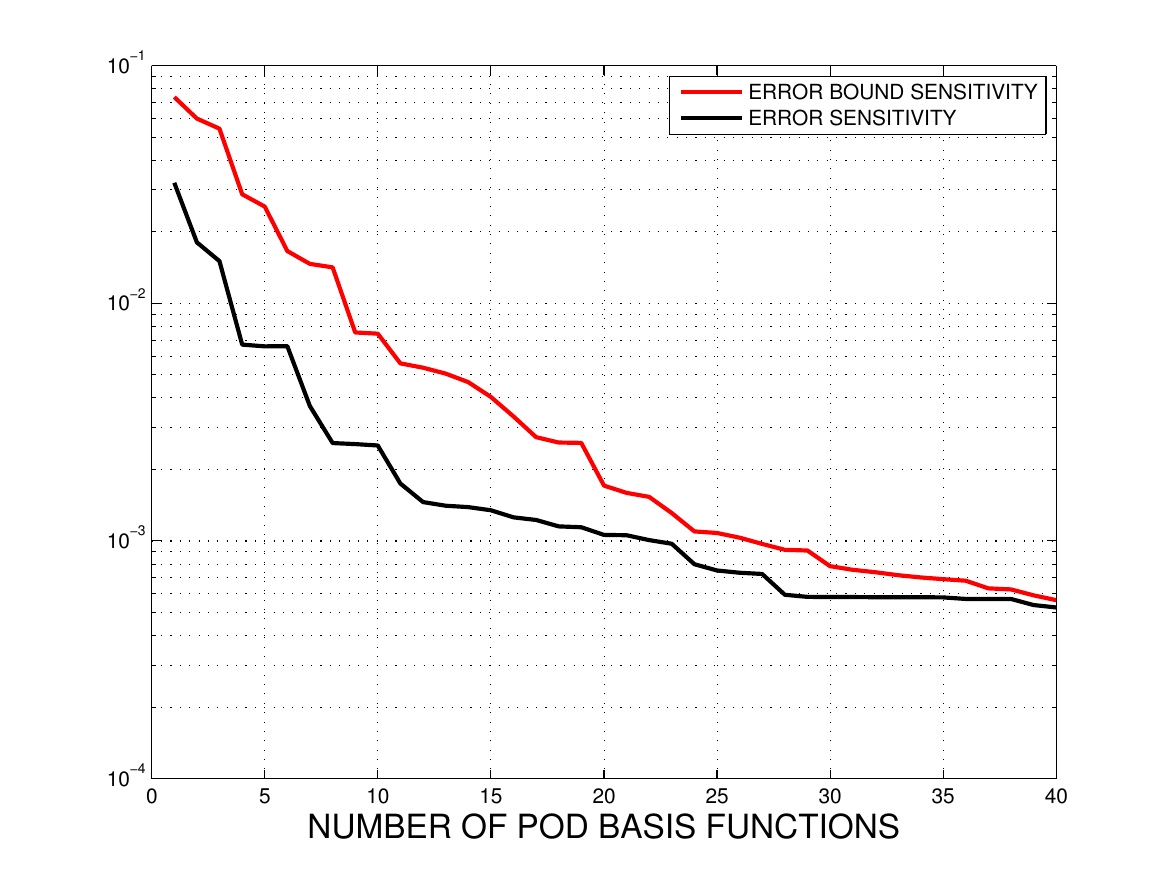}\\
\includegraphics[width=6cm]{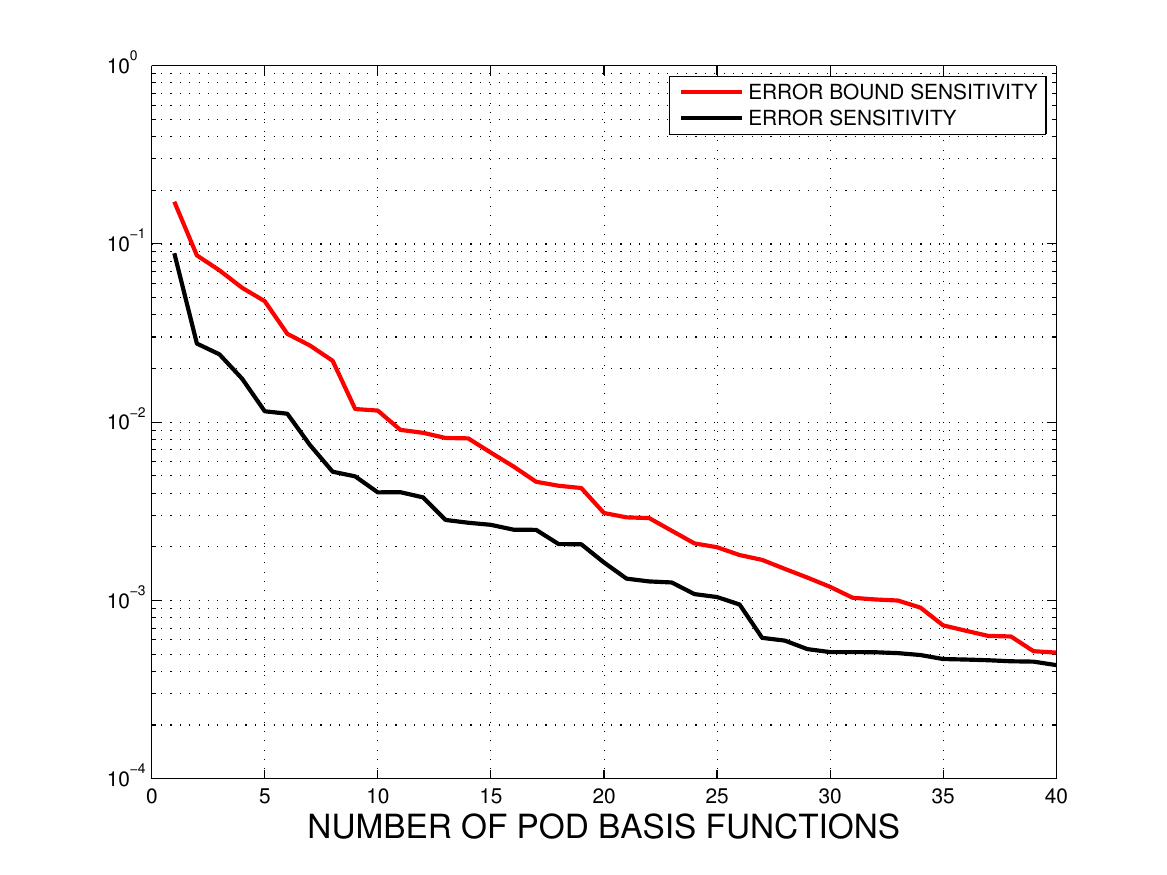}
\includegraphics[width=6cm]{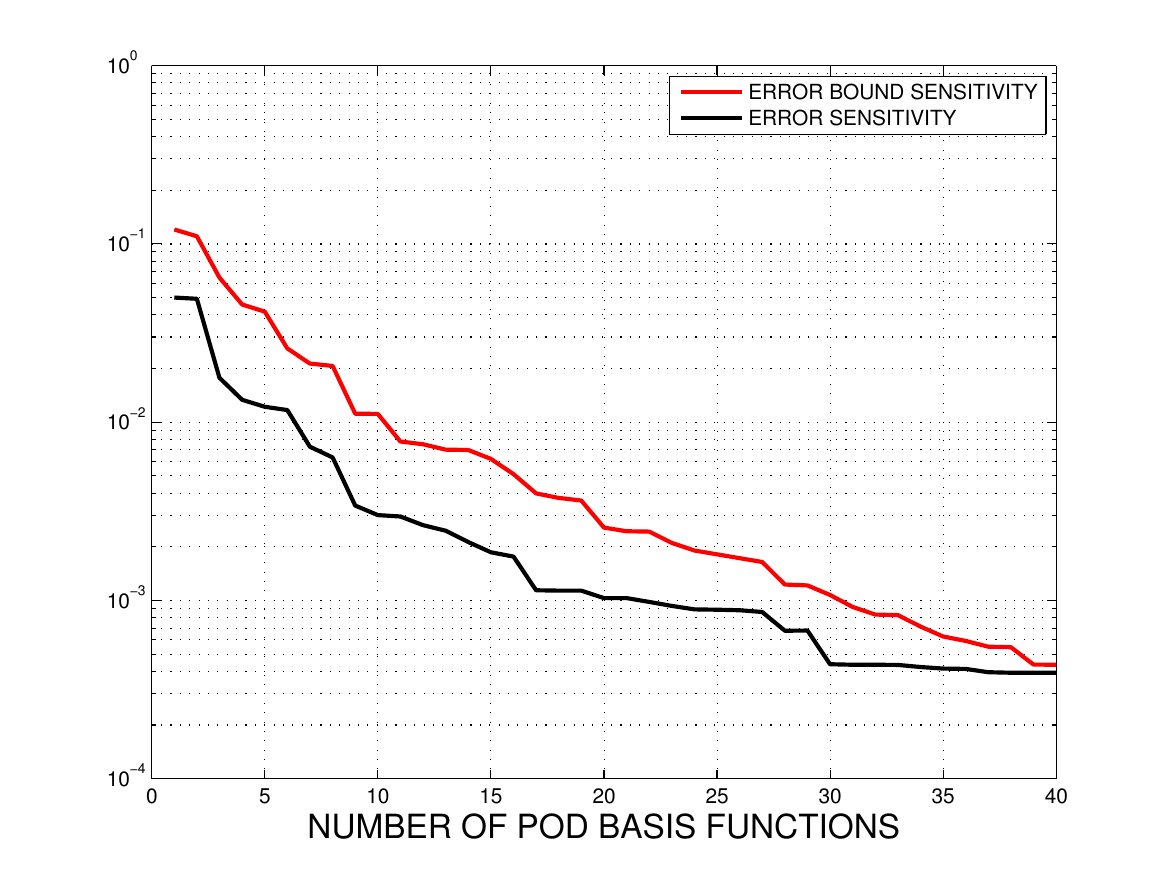}
\caption{\label{fig:apost} Maximum error over all parameter configurations related to $\mathcal{D}_{test}$.  Error behavior $\mathcal{E}_0$ and error bound for state equation (top-left), Error behavior $\mathcal{E}_{|\mu|}$ for $|\mu|=1$ and error bound for the sensitivity equations (top-right and bottom).}
\end{figure}

In Figure \ref{fig:apost}, we present an error analysis to check the quality of our a-posteriori error discussed in Section \ref{sec:pmor}. We compare the error between the POD solution and the high dimensional approximation for state and sensitivities with the error bound presented in Section 4.2. 
With the {\bf W}-norm introduced in Section \ref{sec:pmor}, we set
$$\mathcal{E}_{|\mu|}(p):=\max_{p\in\mathcal{D}_{test}} \|  {\bf u}^n_{h,\mu}(p)  - {\bf u}^{n,\ell}_{h,\mu}(p)\|_{\bf W},\qquad |\mu| =0,1,$$

%$$\mathcal{E}_n(p):=\max_{p\in\mathcal{D}_{test}} \|  {\bf u}^n_{h,\mu}(p)    \bf{u}_h^n(p)-\bf{u}_h^{n, \ell}(p)\|_{\bf W},\qquad n=0,1,2,3.$$

where $\mathcal{D}_{test}$ is chosen in the centers of the boxes of $\mathcal{D}_{train}$ as follows:
$$ \mathcal{D}_{test} = \{3.375, 8.125, 12.875, 17.625\} \times \{1.5, 2.5, 3.5, 4.5\} \times \{4.75,  6.25, 7.750, 9.25\}.$$ 
As one can see model order reduction is able to reach only an error of size $10^{-4}$ (with a large amount of basis functions) and we would need to increase the number of basis functions to achieve higher accuracy.
We also want to emphasize that a standard POD approach, with only state snapshots in $\mathcal{D}_{train}$, is able to reach an accuracy of order $10^{-2}$ with $50$ basis functions.

\begin{figure}[htbp]
\centering
\includegraphics[width=6cm]{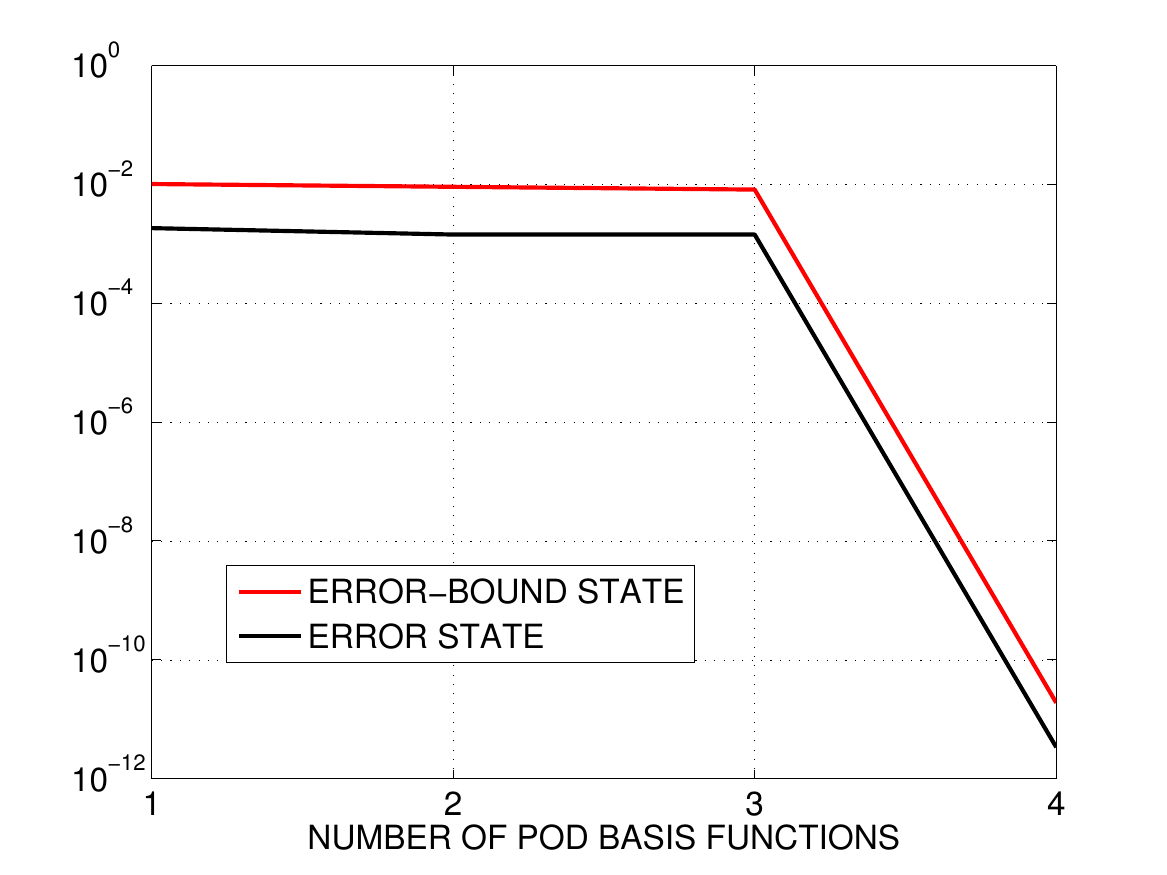}
\includegraphics[width=6cm]{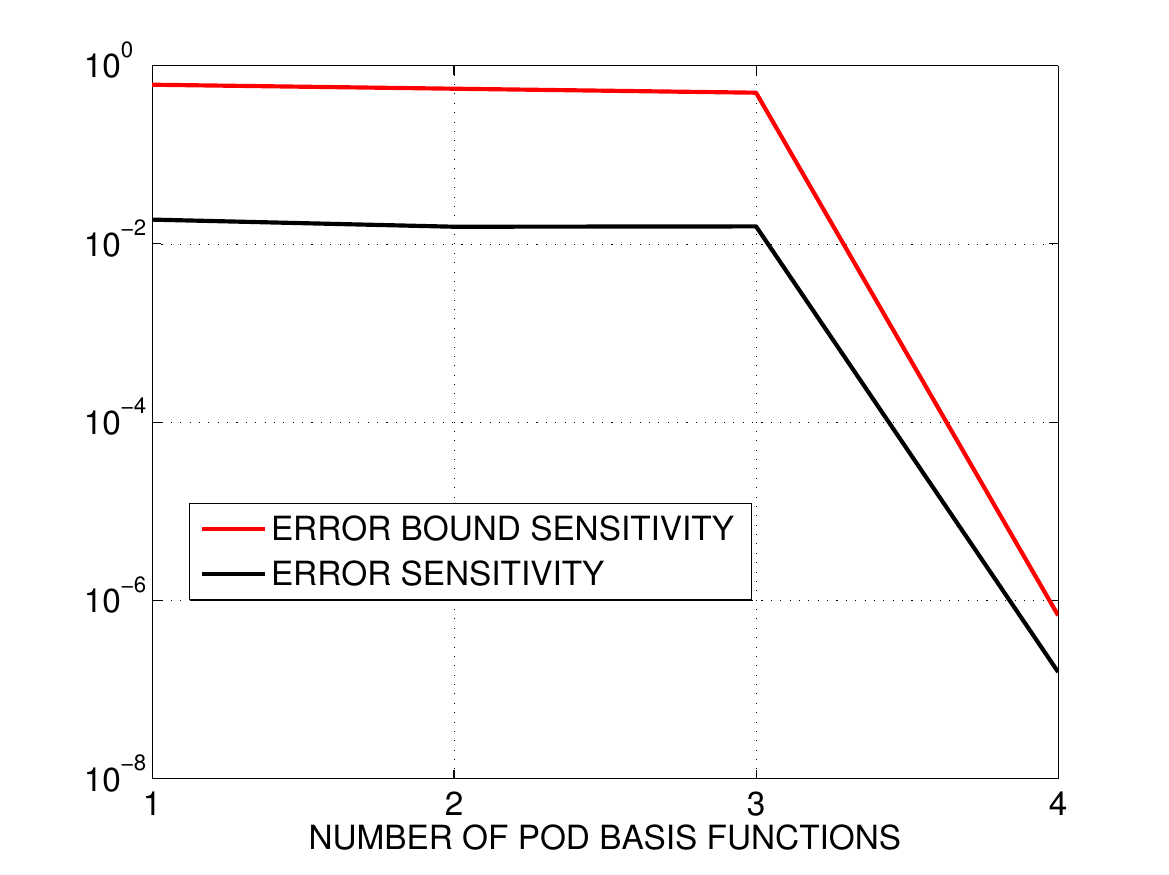}\\
\includegraphics[width=6cm]{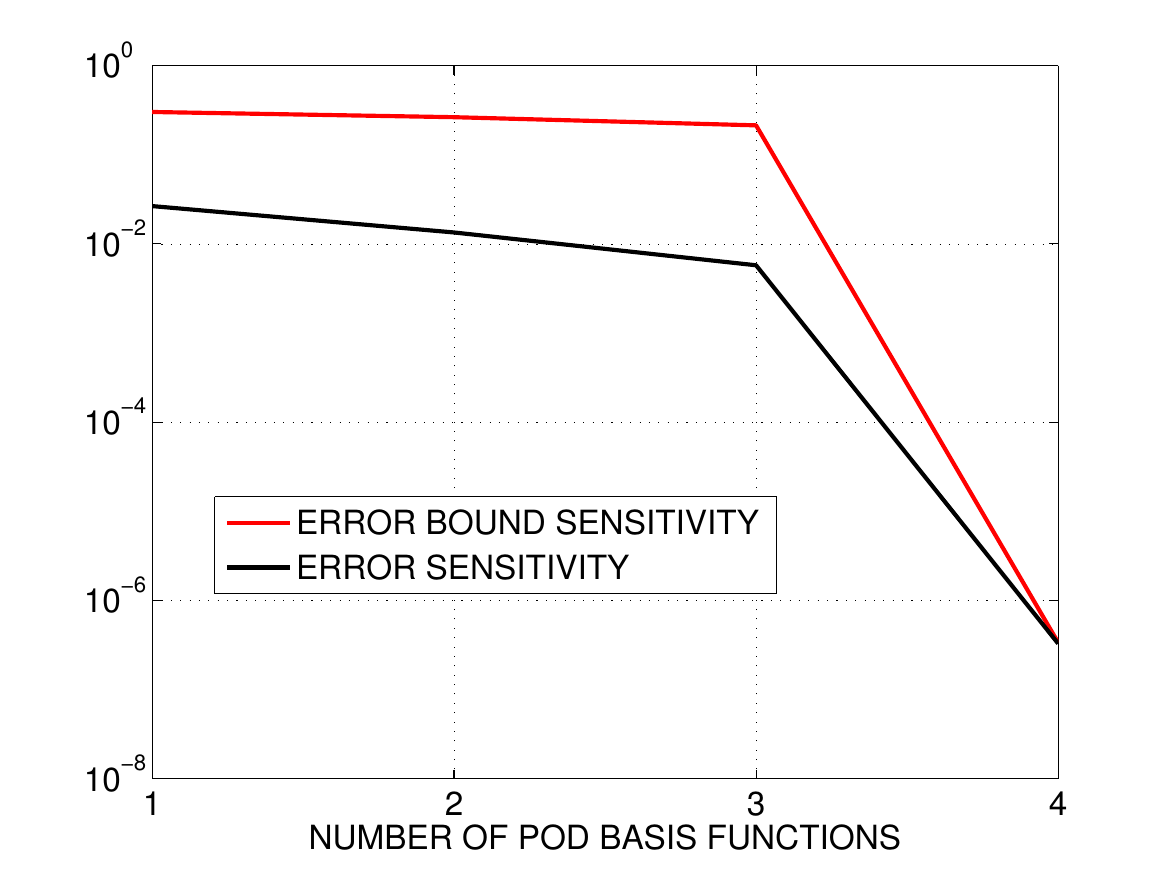}
\includegraphics[width=6cm]{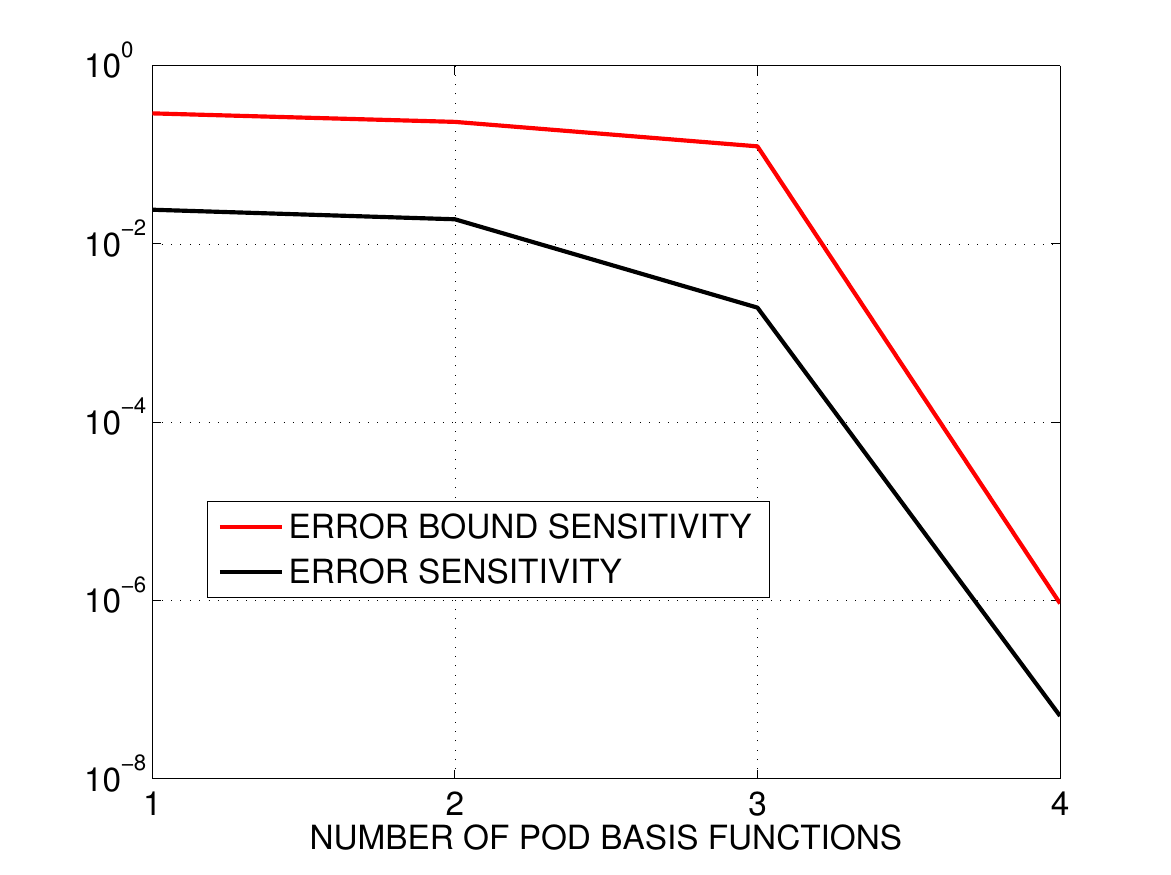}
\caption{\label{fig:apost2}  {Maximum error of the surrogate model computed with respect to a parameter in the neighborhood of $\mathcal{D}_{train}$.  Error behavior for state equation  (top-left), and for the sensitivity equations for each parameter (top-right and bottom).}}
\end{figure}

In Figure \ref{fig:apost2}, we show \blue{that} it is possible to reach an accuracy of order $10^{-6}$ or higher with only 4 POD basis functions if we compute the snapshots with respect to one parameter in $\mathcal{D}_{train}$ and then compute the error in a ball of radius 0.1, centered in the chosen parameter. This is our motivation to introduce \blue{Algorithm} 1 for the successive enrichment of the POD model. In fact, we start our algorithm with only one parameter and we require the combined snapshots since in any other case we will not have enough data to generate a surrogate model from the simulation and the sensitivities. We note that if we add the sensitivities, the POD basis functions are improved and the parameter domain can be better explored. We refer the interested reader to \cite{HBP09} for more details.

Let us now draw our attention to the optimization problem and its performances when combined with model order reduction. Table \ref{tablenompod} shows the convergence of Algorithm \ref{Alg:POD} for the nominal optimization. The number of updates of the snapshot set is given in the first column. Since our goal is to reach a desired electromotive force, we will also compare $E_0(p,{\bf u}_{h,\mu}^{0,\ell}(p),90)$ obtained from the POD, and $E_0(p,{\bf u}_{h, \mu}^0(p),90)$ obtained from the full simulation. %Furthermore, we note that the electromotive force is a linear function which allows us to use the a-posteriori estimator for the state and \new{sensitivity} variables. 
As one can see with the surrogate-based optimization obtain the same results as with the full model. The fifth column presents the number of iterations needed in each sub-optimization problem and the sixth column presents the worst case of the error estimate $\Delta_{u_\mu^n}^{\ell}(p)$ for $0\leq |\mu| =n\leq 1$. In the stopping criterium of Algorithm \ref{Alg:POD} we choose $tol=1e-4$. In the last column we show the number of basis functions which also represents the dimension of the reduced problem. The dimension of the reduced space corresponds to the rank of current snapshot \blue{matrix}. We note that the results of the algorithm are close to what we have shown in Table \ref{tab:results2}. The CPU time will be discussed at the end of the section.

\begin{table}[htbp]
\centering
\begin{tabular}{cccccccc}
\toprule
iter & $E_0(p,{\bf u}_h(p),90)$ & $E_0(p,{\bf u}_h^\ell(p),90)$ & $V$ & $\# it$ & $\max\limits_{|\mu|\le 1} \Delta_{u_\mu^n}^{\ell}(p^k)$ & $\ell $\\
\midrule
  0  & 30.3483 & 30.3483 & 133.00 & &\\
  1  &30.3764 & 30.3483 & \phantom{x}62.82 & 8 &  3.26e-02 & 4\\
 2 & 30.3483 & 30.3483 & \phantom{x}62.36  & 4 &  2.97e-06 & 8\\
 \bottomrule
   \end{tabular}
\caption{Performance of POD with nominal optimization}
\label{tablenompod}
\end{table}

Next we analyze the performance of POD for the robust optimization problem. In Table \ref{tablelinpod} we present the results of POD with a linear approximation of the robust counterpart with $n_{max}=2$. As one can see after two iterations we are able to recover the same results obtained without model order reduction, obtaining also a good approximation for the electromotive force (compare column 2 and column 3 in Table \ref{tablelinpod}).

\begin{table}[htbp]
\begin{center}
\begin{tabular}{cccccccc}
\toprule
iter & $E_0(p,{\bf u}_h(p),90)$ & $E_0(p,{\bf u}_h^\ell(p),90)$ & $V$ & $\# it$ &$\max\limits_{|\mu|\le 2} \Delta_{u_\mu^n}^{\ell}(p^k)$ & $\ell $\\
\midrule
  0  & 30.3483 & 30.3483 & 133.00 & &\\
  1  & 30.7319 & 30.5729& \phantom{x}69.41 &8 & 7.12e-01 & 8\\
 2 & 30.5817 & 30.5816 & \phantom{x}66.46 & 4 & 9.83e-05 &16\\
 \bottomrule
   \end{tabular}
\end{center}
\caption{ Performance of POD with linear robust optimization}
\label{tablelinpod}
\end{table}

\blue{Then, we present the results of the surrogate models} for the solution of the robust optimization problem with quadratic approximation with $n_{max}=3$ in Table \ref{tablequadpod} and with \new{quadratic approximation with moving expansion point} in Table \ref{tablequadshiftpod}. The same considerations discussed in the linear example hold true. We note that the number of POD basis functions used in each optimization in this case is larger due to a richer snapshot set.

\begin{table}[htbp]
\centering
\begin{tabular}{cccccccc}
\toprule
iter & $E_0(p,{\bf u}_h(p),90)$ & $E_0(p,{\bf u}_h^\ell(p),90)$ & $V$ & $\# it$ & $\max\limits_{|\mu|\le 3} \Delta_{u_\mu^n}^{\ell}(p^k)$ & $\ell $\\
\midrule
  0  & 30.3483 & 30.3483 & 133.00 & &\\
  1  & 30.8112 & 30.7332 & \phantom{x}68.12 &6 & 3.32e-02 & 8\\
  2 & 30.6997 & 30.6995 & \phantom{x}67.96 & 5 & 8.07e-05& 16\\
  \bottomrule
   \end{tabular}
\caption{Performance of POD with quadratic robust optimization}
\label{tablequadpod}
\end{table}

\begin{table}[htbp]
\centering
\begin{tabular}{cccccccc}
\toprule
iter & $E_0(p,{\bf u}_h(p),90)$ & $E_0(p,{\bf u}_h^\ell(p),90)$ & $V$ & $\# it$ & $\max\limits_{|\mu|\le 3} \Delta_{u_\mu^n}^{\ell}(p^k)$ & $\ell $\\
  \midrule
  0  & 30.3483 & 30.3483 & 133.00 & &\\

1  & 30.7882 & 30.7384 &  \phantom{x}67.59 & 9 & 4.05e-02 & 8\\
  2  & 30.7908 & 30.7913 &  \phantom{x}67.74 & 9 & 7.11e-04 & 15\\
   3  & 30.7934 & 30.7933 &  \phantom{x}67.87 & 9 & 7.29e-05 &16 \\
   \bottomrule
   \end{tabular}
\caption{\new{Performance of POD with quadratic robust optimization with shift}}
\label{tablequadshiftpod}
\end{table}

Also in the quadratic case we obtain the same approximation quality as in the full model at lower computational costs (\blue{see Table \ref{tablequadpod}} \new{and Table \ref{tablequadshiftpod}}).

\new{To summarize, all the performances of the surrogate models are presented in Table~\ref{tab:results_new}. Again, we emphasize the results agree with Table~\ref{tab:results1}.} %In particular we note the worst-case approximation with POD better match the results presented without model order reduction and shows}. %Furthermore, the presented numbers for both the linear and the quadratic case are in very good agreement with those presented in Table \ref{tab:results1}.
% ATTENTION: Here we need a transition.
\new{Strictly speaking, none of the solutions obtained with reduced models are robust for the high-fidelity problem.
This is demonstrated by the last column of Table~\ref{tab:results_new}.
All of the numbers should be, but none of them is greater than or equal to $30.34$.
The relative worst-case infeasibility violation ranges from $1.9\,\%$ for the nominal solution to less than $0.07\,\%$ for the improved quadratic model.
Seeing that the reduced models are able to reduce the volume more than the high-fidelity models, this is to be expected.
In many applications, an infeasibility as low as $0.07\,\%$ in the worst case presumably is acceptable; if not, a post-optimization analysis is inevitable.}

%Even in this case we are able to have the same approximation of the full problem in a less expensive framework from the computational point of view. As we can see the results in Table \ref{tab:results1} match perfectly both the linear and quadratic approximation; even the worst-case computed from the surrogate model leads the same results.

\begin{table}
\centering
\begin{tabular}{lccccc}
 \toprule
            & $V_{\mathrm{pm}}$ & \%  & $p$ & $E_0$ ($90$) & $E_0^{\mbox{worst}}$\\
 \midrule
Init    & 133.00 & 100 & (19.00  7.00  7.00) & 30.34 & 29.88\\
 Nom. Opt  &  61.24 &  46\%   & (21.25  2.88  6.87) & 30.23 & 29.77\\
 Rob. Lin     &  66.71 &  50\%   & (20.76  3.21  6.14) & 30.56 & 30.21\\
 Rob. Quad & 67.48  & 51\%   & (21.21  3.17  6.92) & 30.76  & 30.29\\
 Rob. Shift  & 67.86  & 51\%   & (21.31  3.18  6.96) & 30.79  & 30.32\\
  \bottomrule
\end{tabular}
\caption{\new{Comparison of the results obtained by POD coupled with the optimization and robust optimization using different approximation orders.}}
\label{tab:results_new}
\end{table}

%\myc{Check CPU and it}

The computational costs of the POD approach are summarized in Table \ref{tab:podresults}. As one can see we obtain \red{a reduction} of the amount of PDEs solved (compare Table \ref{tab:results2}), keeping the same accuracy attained in the high dimensional problem. For instance, in the quadratic robust optimization with model order reduction we only have to perform $36$ PDE solves in contrast to $72$ PDE solves in the optimization with the high-fidelity finite element model. \new{The POD approach for quadratic approximation with shift presents a speed up of factor $5$.}
Although we solve more (reduced) PDEs the CPU time is reduced due to the number of full dimensional PDEs to be solved in the full setting, compare Tables \ref{tab:results2} and \ref{tab:podresults}.  Finally, we note that the CPU time involves both offline and online stage, so that the speed up is considerable.

\begin{table}[htbp]
\centering
\begin{tabular}{lcccc}
 \toprule
             & CPU (s)  & PDE solves & Reduced PDE solves\\
 \midrule
% Nom. Opt   & $ 52.29 $  & $9$  & $45$\\
Nom. Opt & $7.7$ &  $8$ & $48$  \\
 Rob. Lin   & $11.3$  &  $24$ &$110$\\
 Rob. Quad  & $14.2$   & $36$ & $132$\\
 \new{Rob. Shift}& $20.9$  &  $48$ & 324\\
 \bottomrule
\end{tabular}
\caption{Performance of the POD method and computational cost.}
\label{tab:podresults}
\end{table}
 
% -------------------------------------------------
% -------------------------------------------------
% ------------- CONCLUSIONS --------------
% -------------------------------------------------
% -------------------------------------------------

\section{Conclusion}
In this paper we present a new approach which combines model order reduction to robust optimal control. 
We investigate a \red{parameter optimization problem} governed by a parametric elliptic partial
differential equation with uncertain parameters. We introduce a robust optimization framework that
accounts for uncertain model and optimization parameters.
\new{By shifting the expansion point of the utilized Taylor approximation, we are able to model the nonlinear effects of the uncertain parameter more accurately and improve the robust solutions slightly.}
The resulting non-linear optimization problem has a bi-level structure due to the min-max formulation.
%The resulting optimization problem,
%then, has a bi-level structure for the solution of this problem which leads to a non-linear optimization
%problem with a min-max formulation. 
We propose an adaptive model order reduction approach which avoids long offline stages and provides a certified
reduced order surrogate model for the parametrized PDE which then is utilized in the numerical
optimization. The presented numerical results clearly illustrate the validity and performance of the presented approach.
%\myc{add comments about shift}
\section{Acknowledgements}
The authors wish to acknowledge the support of the German BMBF in the context of the SIMUROM project (grant number 05M2013)\new{, and the support of the German Research Foundation in the context of SFB~805}.
%Fist author wishes also to acknowledge the support of the Department of Energy grant number DE-SC0009324.


\begin{thebibliography}{30}

\bibitem{AH01}
K. Afanasiev and M. Hinze. 
\newblock {\em Adaptive control of a wake flow using proper orthogonal decomposition.}
\newblock Lecture Notes in Pure and Applied Mathematics, {\bf 216}, 2001, 317--332. 
%\newblock Shape Optimization \& Optimal Design, Marcel Dekker, 2001. 

\bibitem{AHLU15}
A. Alla, M. Hinze, O. Lass and S. Ulbrich.
\newblock {\em Model order reduction approaches for the optimal design of permanent magnets in electro-magnetic machines.}
\newblock IFAC-PapersOnLine, {\bf 43}, 2015, 242-247.

\bibitem{AU16}
A. Alla, U. Matthes.
\newblock{\em Model order reduction for a linearized robust PDE constrained optimization.}
\newblock IFAC-PapersOnLine, {\bf 49}, 2016, 321-326.

%\bibitem{AV15}
%A. Alla and S. Volkwein.
%\newblock {\em Asymptotic stability of POD based model predictive control for a semilinear parabolic PDE.}
%\newblock  Adv. Comput. Math. {\bf 41}, 2015, 1073-1102. 

\bibitem{A05} 
A.C. Antoulas.
\newblock Approximation of Large-Scale Dynamical Systems.
\newblock SIAM, 2005.

\bibitem{AFS02}
E. Arian, M. Fahl and E. Sachs.
\newblock {\em Trust-region proper orthogonal decomposition models by optimization methods.}
\newblock In Proceedings of the 41st IEEE Conference on Decision and Control, Las Vegas, Nevada, 2002, 3300--3305.

%\bibitem{AK01}
%J.A. Atwell and B.B. King. {\em Proper orthogonal decomposition for reduced basis feedback controllers for parabolic equations}, Matematical and Computer Modelling, \textbf{33}, 2001, 1-19. 

\bibitem{BGN09}
A.~Ben-Tal, A.~Goryashko and A.~Nemirovski.
\newblock \emph{Robust optimization}.
\newblock Princton University Press, 2009.

\bibitem{BN02}
A.~Ben-Tal and A.~Nemirovski. 
\newblock {\em Robust optimization -- methodology and application.} 
\newblock Mathematical Programming, {\bf 92}, 2002, 453--480.

\bibitem{BBC11}
D.~Bertsimas, D.B.~Brown and C.~Caramanis.
\newblock {\em Theory and application of robust optimization.}
\newblock SIAM Review, {\bf 53}, 2011, 464--501.

\bibitem{BL97}
J.~Birge and F.~Louveaux.
\newblock {\em Introduction to Stochastic Programming}.
\newblock Springer, 1997.

\bibitem{CGT00}
A.R.~Conn, N.I.M.~Gould and P.L.~Toint.
\newblock {\em Trust-Region Methods}. 
\newblock Society for Industrial and Applies Mathematics, 2000.

\bibitem{DBK06}
M.~Diehl, H.G.~Bock and E.~Kostina. 
\newblock {\em An approximation technique for robust nonlinear optimization.} 
\newblock Mathematical Programming, {\bf 107}, 2006, 213--230.

\bibitem{DGMM08}
M.~Diehl, J.~Gerhard, W.~Marquardt and M.~M\"onnigmann. 
\newblock {\em Numerical solution approaches for robust nonlinear optimal control problems.} 
\newblock Computer \& Chemical Engineering, {\bf 32}, 2008, 1287--1300.

\bibitem{DH15}
  M.~Dihlmann and B.~Haasdonk.
  \newblock {\em Certified {PDE}-constrained parameter optimization using reduced basis surrogate models for evolution problems.}
  \newblock Computational Optimization and Applications, {\bf 60}, 2015, 753--787.

\bibitem{Eva08}
 L.C. Evans.
\newblock {\em Partial Differential Equations}. 
\newblock American Math. Society, Providence, Rhode Island, 2008.

\bibitem{FLRS06}
R.~Fletcher, S.~Leyffer, D.~Ralph and S.~Scholtes.
\newblock {\em Local convergance of SQP methods for mathematical programs with equilibrium constraints.}
\newblock SIAM Journal on Optimization, {\bf 17}, 2006, 259--286.

\bibitem{GV13}
M. Gubisch and S. Volkwein.
\newblock {\em Proper Orthogonal Decomposition for Linear-Quadratic Optimal Control.}
\newblock In P. Benner, A. Cohen, M. Ohlberger, and K. Willcox (eds.), Model Reduction and Approximation: Theory and Algorithms. 5-64, SIAM, Philadelphia, PA, 2017.

\bibitem{H17}
B. Haasdonk.
\newblock Reduced basis methods for parametrized PDEs - a tutorial introduction for stationary and instationary problems.
\newblock  Chapter in P. Benner, A. Cohen, M. Ohlberger and K. Willcox (eds.): Model Reduction and Approximation: Theory and Algorithms, \blue{65-136}, SIAM, Philadelphia, 2017.

\bibitem{HO08}
B. Haasdonk, and M. Ohlberger.
\newblock Reduced Basis Method for Finite Volume Approximations of Parametrized Linear Evolution Equations.
\newblock M2AN, Math. Model. Numer. Anal.,{\bf 42}, 2008, 277--302.

\bibitem{HBP09}
A. Hay, J. Borggaard, D. Pelletier.
\newblock {\em Local improvements to reduced-order models using sensitivity analysis of the proper orthogonal decomposition,} J. Fluid Mech. {\bf 629}, 2009, 41--72. 

\bibitem{HPUU09} 
M. Hinze, R. Pinnau, M. Ulbrich and S. Ulbrich. 
\newblock {\em Optimization with PDE Constraints. Mathematical Modelling: Theory and Applications, 23}. 
\newblock Springer Verlag, 2009.

%\bibitem{HV05}
%M. Hinze and S. Volkwein. {\em Proper orthogonal decomposition surrogate models for nonlinear dynamical systems: error estimates and suboptimal control}, in Reduction of Large-Scale Systems, P. Benner, V. Mehrmann, D. C. Sorensen (eds.), Lecture Notes in Computational Science and Engineering, \textbf{45}, 2005, 261-306.

\bibitem{HPHB97}
S. Hennerberger, U. Pahner, K. Hameyer and R. ~Belmans.
\newblock {\em Computation of a highliy satured permanent magnet synchronous motor for a hybrid electric vehicle.}
\newblock IEEE Trans. Magn., {\bf 33}, 1997, 4086--4088.

\bibitem{HD13}
B.~Houska and M.~Diehl.
\newblock {\em Nonlinear robust optimization via sequential convex bilevel programming.}
\newblock Mathematical Programming, {\bf 142}, 2013, 539--577.

\bibitem{KTV13}
E. Kammann, F.~Tr\"oltzsch and S. Volkwein. 
\newblock {\em A method of a-posteriori error estimation with application to proper orthogonal decomposition.}
\newblock ESAIM: M2AN, {\bf 47}, 2013, 555--581.

% \bibitem{KLU17}
% P. Kolvenbach, O. Lass and S. Ulbrich. 
% \newblock An approach for robust PDE-constrained optimization with application to shape optimization of electrical engines and of dynamic elastic structures under uncertainty.
% \newblock {\em Submitted} 2017.

\bibitem{KG14} 
\blue{M. K\"archer} and M. Grepl.
\newblock {\em A certified reduced basis method for parametrized elliptic optimal control problems.}
\newblock ESAIM: COCV, {\bf 20}, 416--441, 2014.

\bibitem{KS16}
D.P.~Kouri and T.M.~Surowiec.
\newblock {\em Risk-averse PDE-constrained optimization using the conditional value-at-risk.}
\newblock SIAM Journal on Optimization, {\bf 26}, 2016, 365--396.

%\bibitem{KVX04}
%K. Kunisch, S. Volkwein, and L. Xie. 
%\newblock {\em HJB-POD based feedback design for the optimal control of evolution problems}, 
%\newblock SIAM J. on Applied Dynamical Systems, \textbf{4}, 2004, 701-722.

\bibitem{LU17}
O. Lass and S. Ulbrich.
\newblock {\em Model order reduction techniques  with a posteriori error control for nonlinear robust optimization governed by partial differential equations.}
\newblock SIAM J. on Sc. Comp., {\bf 39}, 2017, S112-S139.

\bibitem{Ley06} 
S. Leyffer. {\em Complementarity constraints as nonlinear equations: Theory and numerical experience},
 in Optimization with Multivalued Mappings: Theory, Applications, and \blue{Algorithms}, S. Dempe and V. Kalashnikov (eds.), {\bf 2} of Springer Series in Optimization and
Its Applications, Springer, 2006, 169--208.

\bibitem{NW06} 
J. Nocedal and S.J. Wright. 
\newblock \emph{Numerical Optimization, second edition}. 
\newblock Springer Series in Operation Research and Financial Engineering, 2006.

\bibitem{OH13}
P.~Offermann and K.~Hameyer.
\newblock {\em A polynomial chaos meta-model for nonlinear stochastic magnet variations.}
\newblock COMPEL - The international journal for computation and mathematics in electrical and electronic engineering, {\bf 32}, 2013, 1211--1218.

\bibitem{Oks05}
B.~\O{}ksendal.
\newblock {\em Optimal control of stochastic partial differential equations.}
\newblock Stochastic Analysis and Applications, {\bf 23}, 2005, 165--179.


\bibitem{OP07} I.B. Oliveira and A.T. Patera.
\newblock{\em Reduced-basis techniques for rapid reliable optimization of systems described by affinely parametrized coercive elliptic partial differential equations.}
\newblock Optim Eng {\bf 8}, 2007, 43--65.


\bibitem{P98} 
U. Pahner.
\newblock {\em A general design tool for theorical optimization of electromagnetic energy transducers.} 
\newblock PhD Thesis, KU Leuven, 1998.

\bibitem{PR07}
A. T. Patera and G. Rozza.
\newblock {\em Reduced Basis Approximation and A Posteriori Error Estimation for Parametrized Partial Differential Equations}
\newblock MIT book, 2007.

\bibitem{P83}
M.J.D. Powell.
\newblock {\em Variable Metric Methods for Constrained Optimization.}
\newblock {In: Mathematical Programming -- The State of the Art: Bonn 1982, A. Bachem, B. Korte, and M. Grötschel (eds.), Springer, 1983, 288--311}

\bibitem{QGVW16} 
E.~Qian, M.~Grepl, K.~Veroy and K.~Willcox.
\newblock {\em A Certified Trust Region Reduced Basis Approach to PDE-Constrained Optimization.}
\newblock SIAM J. Sci. Comput., {\bf 39}, 2016, S434-S460. 


\bibitem{RZ91} 
M.A. Rahman and P. Zhou.
\newblock {\em Determination of saturated parameters of PM motors using loading magnetic fields.}
\newblock IEEE Trans. Magn., {\bf 27}, 1991, 3947--3950.

\bibitem{RHP08}
G.~Rozza, D.B.P.~Huynh and A.T.~Patera.
\newblock {\em Reduced Basis Approximation and a Posteriori Error Estimation for Affinely Parametrized Elliptic Coercive Partial Differential Equations.}
\newblock Arch. Comput. \blue{Methods} Eng., {\bf 15}, 2008, 229--275.

%\bibitem{SV10}
%E.W. Sachs and S. Volkwein. {\em POD Galerkin approximations in PDE-constrained optimization}, GAMM Mitteilungen, \textbf{33}, 2010, 194-208.

\bibitem{Sic13}
A.~Sichau.
\newblock {\em Robust Nonlinear Programming with Discretized PDE Constraints using Second-order Approximations.}
\newblock PhD Thesis, TU Darmstadt, 2013.

\bibitem{Tro10}
F.~Tr\"oltzsch.
\newblock \emph{Optimal Control of Partial Differential Equations: Theory, Methods and Application},
\newblock American Mathematical Society, 2010.

\bibitem{TKXP12}
H.~Tiesler, R.M.~Kirby, D.~Xiu, and T.~Preusser.
\newblock {\em Stochastic Collocation for Optimal Control Problems with Stochastic PDE Constraints.}
\newblock J. Control Optim., {\bf 50}, 2012, 2659--2682. 

\bibitem{TV09} 
F.~Tr\"oltzsch and S. Volkwein.
\newblock {\em POD a-posteriori error estimates for linear-quadratic optimal control problems.}
\newblock \blue{Comput. Optim. and Appl.}, {\bf 44}, 2009, 83--115.

\bibitem{V13}
S. Volkwein. {\em Model Reduction using Proper Orthogonal Decomposition}, Lecture Notes, University of Konstanz, 2013.

\bibitem{ZF14}
M. Zahr and C. Farhat.
\newblock {\em Progressive Construction of a Parametric Reduced-Order Model for PDE-Constrained Optimization.}
\newblock  Internat. J. Numer. Methods Engrg. {\bf 102}, 2015, 1111--1135. 

\bibitem{Zha07}
Y.~Zhang.
\newblock {\em General robust-optimization formulation for nonlinear programming.}
\newblock \blue{J. Optim. Theory and Appl.}, {\bf 132}, 2007, 111--124.

\end{thebibliography}
\end{document}